\documentclass{article}

\usepackage{graphicx}
\usepackage{amsmath}
\usepackage{amssymb}
\usepackage{inputenc}
\usepackage{xcolor}

\numberwithin{equation}{section}

\newtheorem{definition}{Definition}
\newtheorem{lema}{Lemma}

\newtheorem{thm}{Theorem}

\newcommand{\red}[1]{{#1}}
\newcommand{\blue}[1]{{#1}}

\usepackage{geometry}

\title{A structure preserving numerical scheme for Fokker-Planck equations of structured neural networks with learning rules}
\author{Qing He\footnote{Department of Mathematics, Southern Methodist University, Dallas, TX 75275, USA (andrewho@smu.edu).}, \ Jingwei Hu\footnote{Department of Applied Mathematics, University of Washington, Seattle, WA 98195 (hujw@uw.edu).}, \ and \ Zhennan Zhou\footnote{Beijing International Center for Mathematical Research, Peking University, Beijing, China (zhennan@bicmr.pku.edu.cn).}}

\begin{document}
\maketitle

% \medskip
\begin{abstract}
In this work, we are concerned with a Fokker-Planck equation related to the nonlinear noisy leaky integrate-and-fire model for biological neural networks which are structured by the synaptic weights and equipped with the Hebbian learning rule. The equation contains a small parameter $\varepsilon$ separating the time scales of learning and reacting behavior of the neural system, and an asymptotic limit model can be derived by letting $\varepsilon\to 0$, where the microscopic quasi-static states and the macroscopic evolution equation are coupled through the total firing rate. To handle the endowed flux-shift structure and the multi-scale dynamics in a unified framework, we propose a numerical scheme for this equation that is mass conservative, unconditionally positivity preserving, and asymptotic preserving. We provide extensive numerical tests to verify the schemes' properties and carry out a set of numerical experiments  to investigate the model's learning ability, and explore the solution's behavior when the neural network is excitatory.

\end{abstract}

%\tableofcontents

%\jh{Which abstract the final version?}

%\hq{The latter one is modified by Professor Zhou}

%\jh{My change is highlighted in blue. Many minor changes are done without highlighting. -Jingwei}

%\hq{My respose and questions are highlighted in red. Changes that are entirely accepted are not highlighted. -Qing}

\section{Introduction}
In biological neural networks, neurons transmit information in two ways: one is conveying rapid and transient electrical signals on their membranes, and the other is sending neurotransmitters at a synapse. The learning and memory of the network is based on the connection weights of the synapses among the neurons \cite{Synaptic_plasticity_and_memory}. How a neural network with a given structure reacts to the environment and how the neural network's infrastructure is modified by the environment in the long run are two of the essential topics in neural science.

Mathematically, there has been plenty of models describing the neural networks' behavior from all scales and sizes. The most successful single-neuron model is undoubtedly the Hodgkin-Huxley model \cite{Hodgkin1952a} and its simplification, the Fitzhugh-Nagumo model \cite{FitzHugh1961Impulses}. However, when modeling large-scale biological neural networks, these models are far too complicated and impractical. Therefore, scientists and mathematicians usually adopt the integrate-and-fire model \cite{lapicque1907recherches} when modeling the system that consists of a large number of neurons \cite{caceres2010analysis,carrillo2013classical, Global-in-time,brunel2000dynamics,DGGP17theMeanField,brunel1999fast}. One of the most widely used version of the integrate-and-fire model is the nonlinear noisy leaky integrate and fire model (NNLIF). In this model, when the firing event does not occur, \red{the membrane potential $V(t)$ of a neuron} is influenced by a relaxation to a resting potential $V_{I}$, an incoming synaptic current $\mathcal{I}(t)$ from other neurons and the background. The current $\mathcal{I}(t)$ is approximately decomposed into a drift term and a term of Brownian motion. The SDE, in the simplest form, is given by
\begin{equation}\label{neuron-SDE}
\mathrm d V=-\left(V-V_{I}\right) \mathrm d t+\mu \mathrm d t+\sigma \mathrm d B_{t},
\end{equation}
where the parameters $\mu$ and $\sigma$ are determined by the synaptic current $\mathcal{I}(t)$. The distinguished feature of this model is the incorporation of the firing event: when
a neuron reaches a firing potential $V_{F}$, it discharges itself and the potential jumps to a resetting potential $V_{R}$ immediately \blue{\cite{delarue2015global,delarue2015particle}}. Here, we assume, $V_{R}<V_{F}$, and the firing event is expressed by
\begin{equation}\label{neuron-SDE-firing}
V\left(t^{-}\right)=V_{F}, \quad V\left(t^{+}\right)=V_{R}.
\end{equation}
(1.1) and (1.2) constitute the stochastic process for the NNLIF model, which is a SDE coupled with a renewal process. By It\^o’s formula, one can derive the time evolution of the density function $p(v, t)$ \cite{liu2020rigorous,zhou2021investigating}, which represents the probability of finding a neuron at the voltage $v$ and given time $t$:

\begin{equation}\label{NNLIF-Fokker-Planck}
\begin{cases}
{\displaystyle \frac{\partial p}{\partial t} +\frac{\partial }{\partial v} (hp)-a\frac{\partial ^{2} p}{\partial v^{2}} =N( t) \delta ( v-V_{R}) ,\quad t\in ( 0,\infty ) ,v\in ( -\infty ,V_{F})},\\
{\displaystyle p(v,0)=p_{0} (v),\ \ p(-\infty ,t)=p(V_{F} ,t)=0},
\end{cases}
\end{equation}
where
\begin{equation}
    N(t) = -a\frac{\partial p}{\partial v}(V_F,t)
\end{equation}
stands for the firing rate of the network, $a = {\sigma^2}/{2}$ stands for the amplitude of the noise, \red{$\delta(\cdot)$ stands for the Dirac delta function}, and
\begin{equation} \label{NNLIF-h}
h=h(v,N(t)) = -v + I + w N(t)
\end{equation}
stands for the rate of the rising of the neurons' potential, consisting of a decay term $-v$, an external input function $I$ and a group reaction term $wN(t)$. The coefficient $w$ in (\ref{NNLIF-h}) can be either positive or negative. In the simplest scenario, $w$ is a constant. We say that the neural network is excitatory when $w>0$, and inhibitory when $w < 0$.

When it comes to the learning behavior of the biological neural network, the most seminal work is the Hebbian rule \cite{hebb1949organization} stated as a motto ``neurons that fire together wire together''. A very simple mathematical setup of this rule is as follows: assuming that the strength of weights $w_{i j}$ between two neurons $i$ and $j$ increases when the two neurons have high activity simultaneously. For $M$ neurons in interactions, the classical Hebbian rule relates the weights to the activity $N_{i}$ of the neuron $i$ as

$$
\frac{\mathrm d}{\mathrm d t} w_{i j}=k_{i j} N_{i} N_{j}-w_{i j}, \quad 1 \leq i, j \leq M.
$$

In \cite{BP18}, the authors designed a model that integrates the NNLIF model and the Hebbian rule stated above. They considered a heterogeneous population of homogeneous neural networks structured by their synaptic weights $w \in(-\infty,+\infty)$, negative sign standing for inhibitory neurons and positive sign standing for excitatory neurons. They introduced a learning rule in order to modulate the distribution of synaptic weights $H$ and allow the network to recognize some given input signals $I$ by choosing an appropriate heterogeneous synaptic weight distribution $H$ adapted to the signal $I$. They assumed that the subnetworks interact only via the total firing rate $N$, with synaptic weights described with a single \red{variable} $w$. And they gave the following interpretation: all the subnetworks parametrized by $w$ may modulate their intrinsic synaptic weight $w$ with respect to a function $\Phi$ which depends on the intrinsic activity $N(w)$ of the network parametrized by $w$ and of the total activity of the network $\overline{N}$. Then, the proposed generalization of the Hebbian rule consists in choosing
\begin{equation}
\Phi(N(w), \overline{N})=\overline{N} N(w) K(w), \quad \overline{N}=\int_{-\infty}^{\infty} N(w) d w,
\end{equation}
where $K(\cdot)$ represents the learning strength of the subnetwork with synaptic weight $w$. Adding the above choice of the learning rule, they obtained the following equation
% \jh{
\begin{equation}
\label{5013957753484064}
\displaystyle \frac{\partial p}{\partial t}+\frac{\partial}{\partial v}[(-v+I(w)+w \sigma(\overline{N}(t))) p]+\varepsilon \frac{\partial}{\partial w}[(\Phi-w) p]-a \frac{\partial^{2} p}{\partial v^{2}}=N(w, t) \delta\left(v-V_{R}\right),
\end{equation}%}
where $\varepsilon$ stands for a time scale ratio which takes into account that learning is slower than the typical voltage activity of the network.

If we want to study the learning behavior of the model, we may further perform a time rescaling $t \to t/\varepsilon$ for (\ref{5013957753484064}) to arrive at
\begin{equation} \label{auGi4Tx7zhb}
\left\{\begin{array}{ l }
\begin{aligned}
\frac{\partial p}{\partial t} +\frac{\partial }{\partial w} [(\overline{N}( t) N(w,t)K(w)-w)p] & \\
=\frac{1}{\varepsilon }\Bigl\{a\frac{\partial ^{2} p}{\partial v^{2}} -\frac{\partial }{\partial v} [(-v+I(w) & +w\sigma (\overline{N} (t)))p]+N(w,t)\delta ( v-V_{R})\Bigr\}\\
\text{for }( v,w,t) & \in ( -\infty ,V_{F}) \times ( -\infty ,+\infty ) \times ( 0,T),
\end{aligned}\\[30pt]
\displaystyle N (w,t)=-a\frac{\partial p}{\partial v}( V_{F} -0,w,t) \geq 0,\ \ \overline{N} (t)=\int _{-\infty }^{\infty } N(w,t)\mathrm{d} w,
\end{array}\right. 
\end{equation}
\blue{where $\frac{\partial p}{\partial v}( V_{F} -0,w,t)$ means the left  derivative in the $v$ direction} and the corresponding initial and boundary conditions are
\begin{equation} \label{a52V3vU+wub}
p( V_{F} ,w,t) =0,\ \ p(-\infty ,w,t)=0,\ \ p(v,\pm \infty ,t)=0,\ \ p(v,w,0)=p^{0} (v,w). 
\end{equation}
Here $p(v,w,t)$ at a fixed time $t$ describes the probability density of finding a neuron with synaptic weight $w$ and voltage $v$, henceforth the probability density of finding a single neuron with synaptic weight $w$ is defined as
\begin{equation} \label{aeRiM-jHclb}
H( w,t) =\int _{-\infty }^{V_{F}} p( v,w,t)\mathrm{d} v,
\end{equation}
and the normalization condition of probability ensures that $\int _{-\infty }^{+\infty } H( w,t)\mathrm{d} w=1$ holds for all $ t\geq 0$. $N(w,t)$ describes the firing rate of the neurons with synaptic weight $w$, at time $t$ and $\overline N(t)$ describes the entire network's firing rate at time $t$ summing over all the synaptic weights. $a$ is a positive noise coefficient. The model is a multi-scale transport equation, where the $v$-wise transport represents the integrate-and-fire dynamics of the neural network, and the $w$-wise transport represents the network's learning behavior. 

% \cite{torres_dynamics_2020} also gives a model with learning rule that is even closer to Hebbian's. However, the location in this model is just used as a parameter for location, there is no partial derivatives on location parameter appearing in this equation.

% todo: mode citations on 

\blue{Systematic understanding of the model (\ref{auGi4Tx7zhb}) is largely lacking, although some of its properties were proved in \cite{BP18}}. Particularly, the numerical studies are far from sufficient for understanding the possible uses and limitations of this model. In this work, we further explore this model from a numerical perspective. Specifically, we aim to design a structure-preserving numerical scheme for (\ref{auGi4Tx7zhb}), and further explore this model with extensive numerical experiments to investigate possible solution structures and interpretations.

In our work, we propose a numerical scheme with the positivity preserving and asymptotic preserving properties. There are two key ingredients that make our scheme meet the desired properties. The first one is the Scharfetter-Gummel symmetric reconstruction \cite{JinYanFPLAP11} of the $v$-wise convection and diffusion terms in (\ref{auGi4Tx7zhb}). The second one is that we treat the flux shift term $N\delta(v-V_R)$ in (\ref{auGi4Tx7zhb}) implicitly. 
With the symmetric reconstruction and this implicit treatment, we can prove that implementing our scheme means inverting a so-called $M$-matrix at each time step, and the positivity preserving properties of the scheme is independent of the respective ratio between $\Delta t$ and $\Delta v$ or $\varepsilon$. The only grid ratio requirement is that $\Delta t/\Delta w$ has to satisfy the $w$-wise CFL condition, which is a natural requirement. Furthermore, when letting $\varepsilon\to 0$ without adjusting $\Delta t$, we obtain a numerical scheme consistent with the asymptotic behavior of the model, consisting of a $w$-wise convection equation for $H(w, t)$ defined in (\ref{aeRiM-jHclb}) and a $v$-wise quasi-steady equation for $p(v, w, t)$ with given $H(w, t)$.

\red{Numerical methods for Fokker-Planck type equations are abundant in the literature. In recent years many works have been proposed to preserve the structure of the solution such as mass conservation, positivity, entropy decay, and steady state preserving. This often relies on delicate choice of numerical fluxes and suitable time discretization. Without being exhaustive, we mention a few relevant works. Regarding spatial discretization, there are mainly two kinds of methods: one is based on formulating the equation into the transport form and using upwind fluxes \cite{bessemoulin2012finite}. This is later extended to general nonlinear aggregation and diffusion equations \cite{CCH15}. The other is based on formulating the equation into the diffusion form (i.e., Scharfetter-Gummel form). Depending on how to choose the fluxes, it bifurcates to several variants, for example, the Chang-Cooper scheme \cite{CC70} and its generalization \cite{PZ18}; the symmetrized finite difference scheme \cite{JinYanFPLAP11}, and more recently the high order finite difference scheme \cite{HZ22}. Regarding time discretization, high order explicit or semi-implicit schemes (c.f. \cite{boscarino2016high}) can be used. However, we mention that to get positivity preserving implicit schemes for Fokker-Planck type equations, one is basically limited to first order scheme \cite{liu2018positivity, HH20, BCH20,zhang2021unified}. The only exception we are aware of is the second order exponential method in \cite{HS20} which is applicable to the linear Fokker-Planck equation and hard to generalize to more general cases. For the above reason, we limit ourself to first order method in this paper.
}

%\red{There are some existing numerical work related to the  numerical algorithms for Fokker-Planck equations based on Sharfetter-Gummel type fluxes. \cite{bessemoulin2012finite} presents a numerical scheme for nonlinear degenerate parabolic equations which admit an entropy functional that preserves steady-states and remains valid and second-order accurate in space even in the degenerate case, \cite{pareschi2018structure} presents numerical schemes for nonlinear Fokker-Planck equations that preserve the non negativity of the solution, entropy dissipation and large time behavior,  \cite{boscarino2016high} provides a method to implement Runge-Kutta method in PDE problems with stiff operators, \cite{liu2018positivity} presents a scheme for 2D Keller-Segel equations based on a similar symmetrization reformation we present in this article, and \cite{zhang2021unified} presents a numerical scheme for a parabolic model for multispecies ionic fluids which is second-order accurate in space, first-order in time, and satisfies the positivity preserving and mass conservation properties.}

In the numerical experiment part, we verify the scheme's order of accuracy (first order for $w, t$, second order for $v$) and asymptotic preserving properties, test the model's recognition properties presented in \cite{BP18} and further explore the model's behavior when the solution $p(v,w,t)$ is supported in the region $w > 0$. The numerical results suggest that this model shows great potential in distinguishing a general input by producing distinct and unsymmetrical  responses to orthogonal basis functions. When $w>0$, the solution $p(v,w,t)$ may either develop to a steady state supported in a region $w\in [0,A]$ for some positive $A$, or expand its support rapidly towards $w \to \infty$.

The rest of the paper is organized as follows. In Section \ref{sec: The schemes}, we present our schemes in detail, including the equations and their implementations; in Section \ref{sec:properties}, we give a detailed analysis of our schemes' properties, including conservation, positivity preserving properties and asymptotic preserving properties; and in Section \ref{sec:numerical_results}, we present the numerical results including the verification of the schemes' accuracy, asymptotic preserving properties, numerical exploration of the model's learning and recognition ability and the behavior of the model with positive synaptic weights.

% We gather the equations that describe the asymptotic model again to give a general picture.

% \begin{equation*}
%     \begin{cases}
%         \displaystyle p_{0}( v,w,t)=Q_{\overline{N} [H_{0}( \cdot ,t) ]} (v,w)H_{0} (w,t) & \text{(\ref{aSG4SbuwrU})} \\
%         \displaystyle \frac{\partial H_{0}( w,t)}{\partial t} +\frac{\partial }{\partial w}[(\mathcal{\overline{\mathnormal{N}}}[ H_{0}( \cdot ,t)] N_{\overline{N}[ H_{0}( \cdot ,t)]} H_{0}( w,t) K-w) H_{0}( w,t)] =0 & \text{(\ref{amd9LZngKpb})}
%     \end{cases}
% \end{equation*}
% where $Q_{\overline N}, N_{\overline N}$ and $\overline N[H]$ are defined as 
% \begin{equation*}
%     \begin{cases}
%     \displaystyle \frac{\partial }{\partial v}[ (-v+I(w)+w\sigma (\overline{N} )) Q_{\overline{N}}] -a\frac{\partial ^{2} Q_{\overline{N}}}{\partial v^{2}} = N_{\overline{N}} (w)\delta ( v-V_{R}) ,\\
% \displaystyle Q_{\overline{N}}( V_{F} ,w) =0,\ \ N_{\overline{N}} (w)=-a\frac{\partial Q_{\overline{N}}( V_{F} ,w)}{\partial v} ,\ \ \int _{-\infty }^{V_{F}} Q_{\overline{N}} (v,w)\mathrm{d} v=1,
%     \end{cases}\tag{\ref{aW-HpTlo93}}
% \end{equation*}
% and 
% \begin{equation*}
%     \overline{N}[ H] =\int _{-\infty }^{+\infty } N_{\overline{N}[ H]} (w)H(w)\mathrm{d} w, \tag{\ref{aGNK1JP2NFb}}
% \end{equation*}

\section{The Schemes}\label{sec: The schemes}

Our work aims to provide a numerical scheme for (\ref{auGi4Tx7zhb}). In this section, we will give a detailed description of the schemes we designed.

This section is organized as follows: in Section \ref{subsec:discretization} and Section \ref{subsec:scheme}, we introduce the grid point settings and the numerical scheme for (\ref{auGi4Tx7zhb}), including a $v$-wise semi-implicit ($v$-SI) implementation and a $v$-wise fully-implicit ($v$-FI) implementation; and in Section \ref{subsec: matrix form}, we introduce the matrix form of the $v$-SI and $v$-FI schemes, from which we introduce a useful matrix $\boldsymbol M_j^{\overline{N}}$ to formulate further analysis.

\subsection{Grid point settings and discretizations} \label{subsec:discretization}
Our scheme is a grid-based method, and this subsection gives the basic grid point settings of our scheme. 

We assume that $ p( v,w,t)$ vanishes fast enough as $ v\rightarrow -\infty $ or $ w\rightarrow \pm \infty $ such that we truncate the domain into the bounded region
\begin{equation*}
( v,w) \in [ V_{\min} ,V_{F}] \times [ W_{\min} ,W_{\max}] \times [ 0,T_{\max}],
\end{equation*}
and suppose that $ p( v,w,t)$ is negligible out of this region. Then we divide the intervals $ [ V_{\min} ,V_{F}]$, $[ W_{\min} ,W_{\max}]$, $[ 0,T_{\max}]$ into $ n_{v} ,n_{w} ,n_{t}$ equal sub-intervals with size
\begin{equation} \label{aSC1HLrrxZ}
\Delta v=\frac{V_{F} -V_{\min}}{n_{v}} ,\ \ \Delta w=\frac{W_{\max} -W_{\min}}{n_{w}} ,\ \ \Delta t=\frac{T_{\mathrm{max}}}{n_{t}} ,
\end{equation}
respectively, so that the grid points are assigned as follows\begin{equation} \label{aWZBmXPNc+}
\begin{cases}
v_{i} =V_{\min} +i\Delta v & i=0,1,2,\cdots ,n_{v}\\
w_{j} =W_{\min} +j\Delta w & j=0,1,2,\cdots ,n_{w}\\
t^{m} =m\Delta v & m=0,1,2,\cdots ,n_{t},
\end{cases}
\end{equation}
then $ p_{i,j}^{m}$ represents the numerical approximation of $ p\left( v_{i} ,w_{j} ,t^{m}\right)$. We always assume that there is an index $ r$ satisfying
\begin{equation} \label{a5XOf2DShJb}
V_{R} =v_{r},
\end{equation}
which naturally treats the derivative's discontinuity at $ v=V_{R}$. Actually it can also be seen from (\ref{aSC1HLrrxZ}) and (\ref{aWZBmXPNc+}) that $v_{n_v} = V_F$, so the two points related to flux shift are both aligned with the grid point.
% \jh{Isn't $V_F$ aligned with the grid too?} \hq{Yes, it is. Actually one can derive from (\ref{aSC1HLrrxZ}) and (\ref{aWZBmXPNc+}) that $v_{n_v} = V_F$, so I think we don't need to emphasize that.} \zz{we can add a remark on this to improve the readability.}

For the discretization of $ N = -a\partial p/\partial v|_{v=V_F}$, we approximate the derivative with the first-order finite difference
\begin{equation} \label{a5KFaVLjKZb}
N_{j}^{m} =-a\frac{p_{n_{v} ,j}^{m} -p_{n_{v} -1,j}^{m}}{\Delta v} =a\frac{p_{n_{v} -1,j}^{m}}{\Delta v}
\end{equation}
(here we have applied the boundary condition $ p_{n_{v} ,j}^{m} =0$); for the discretization of $\overline{N} (t)=\int _{-\infty }^{\infty } N(w,t)\mathrm{d} w$, we apply the simplest rectangular rule of numerical integration
\begin{equation} \label{aqZY+PFiRGb}
\overline{N}^{m} =\Delta w\sum _{j=0}^{n_{v}} N_{j}^{m},
\end{equation}
which can also be explained as a trapezoidal rule of numerical integration since the boundary values $ N_{0}^{m}$ and $ N_{n_{w}}^{m}$ are supposed to be negligible (zero); furthermore, we also denote the discrete representation of $H( w,t) =\int _{-\infty }^{V_{F}} p( v,w,t)\mathrm{d} v$ as
\begin{equation} \label{ayvWHg2kk2}
H_{j}^{m} =\Delta v\sum _{i=1}^{n_{v}} p_{i,j}^{m}.
\end{equation}

\subsection{The numerical scheme for the  Fokker-Planck equation (\ref{auGi4Tx7zhb})}
\label{subsec:scheme}

We design a scheme with a finite volume construction: for each $ i=0,1,\cdots ,n_{v} -1$, $ j=0,1,\cdots ,n_{w}$, $ m=0,1,\cdots ,n_{t} -1$, we have the difference equation
\begin{equation} \label{a9OxraVzxdb}
\frac{p_{i,j}^{m+1} -p_{i,j}^{m}}{\Delta t} +\frac{\Phi _{i,j+1/2}^{m} -\Phi _{i,j-1/2}^{m}}{\Delta w} =\frac{1}{\varepsilon }\frac{F_{i+1/2,j}^{m} -F_{i-1/2,j}^{m}}{\Delta v},
\end{equation}
or, in an equivalent splitting form
\begin{equation} \label{aCWU7qDHkqb}
\frac{p_{i,j}^{*m} -p_{i,j}^{m}}{\Delta t} +\frac{\Phi _{i,j+1/2}^{m} -\Phi _{i,j-1/2}^{m}}{\Delta w} =0,
\end{equation}
\begin{equation} \label{ayUOikHTLnb}
\frac{p_{i,j}^{m+1} -p_{i,j}^{*m}}{\Delta t} =\frac{1}{\varepsilon }\frac{F_{i+1/2,j}^{m} -F_{i-1/2,j}^{m}}{\Delta v}.
\end{equation}
This splitting form is introduced only for convenience when discussing the scheme's properties in section \red{\ref{subsec: positivity preserving}}. Our scheme is not based on the idea of operator splitting.

% For the $w$-wise transport, we take the following explicit linearized Riemann solver (Roe solver)
For the $w$-wise transport, we take the following explicit flux construction adapted from Godunov's Method (see, for example, (13.24) in \cite{levenqueNumConserv}): %\hq{I reconsidered my method and think that the scheme I designed (especially when it is written as (\ref{aOpH3A49bX})) is more close to Godunov's original method, so I decided to state my scheme as `an adaptation of Godunov's scheme', instead of `Roe scheme' I always mentioned before}
\begin{equation} \label{aOpH3A49bX}
\Phi _{i,j+1/2}^{m} =
\begin{cases}
\begin{cases}
\min\left\{\Phi _{i,j}^{m} ,\Phi _{i,j+1}^{m}\right\} & p_{i,j}^{m} \leq p_{i,j+1}^{m}\\
\max\left\{\Phi _{i,j}^{m} ,\Phi _{i,j+1}^{m}\right\} & p_{i,j}^{m}  >p_{i,j+1}^{m}
\end{cases} & j=0,\cdots ,n_{w} -1\\
0 & j=-1,n_{w}
\end{cases}
\end{equation}
where
\begin{equation} \label{aCv69yRawS}
\Phi _{i,j}^{m} =(\overline{N}^{m} N_{j}^{m} K( w_{j}) -w_{j} )p_{i,j}^{m} \ \ \ \text{for} \ \ j=0,\cdots ,n_{w} .
\end{equation}

% \jh{I think this is the right place to add CFL condition.}

For the $ v$-wise transport, we need to treat the operator partially implicitly in an appropriate way such that the scheme's stability is not limited by the grid ratio $\Delta t/(\Delta v)^2$ and the stiffness introduced by the smallness of $\varepsilon$. We inherit the idea from \cite{JH19} which implements a Scharfetter-Gummel symmetric reconstruction for the convection-diffusion operator and imposes a flux-shift from the boundary $ V_{F}$ back to $ V_{R}$, but we treat the flux shift operator $N \delta(V-V_R)$ implicitly:
\begin{equation} \label{aiChucap+Ob}
F_{i+1/2,j}^{m} =\begin{cases}
\displaystyle a\frac{M_{i+1/2,j}^{\overline{N}^{\alpha }}}{\Delta v}\left(\frac{p_{i+1,j}^{m+1}}{M_{i+1,j}^{\overline{N}^{\alpha }}} -\frac{p_{i,j}^{m+1}}{M_{i,j}^{\overline{N}^{\alpha }}}\right) +N_{j}^{m+1} \eta ( v_{i+1/2} - V_R) & i=0,1,\cdots ,n_{v} -2,\\
0 & i=-1,n_{v} -1,
\end{cases} \quad 
\end{equation}
where $\eta(\cdot)$ is the Heaviside function ($\eta(x) = 1$ when $x \geq 0$ and $\eta(x) = 0$ when $x < 0$), and
\begin{equation} \label{aqSw71E7wnb}
M_{i,j}^{\overline{N}^{\alpha }} =\mathrm{e}^{-\frac{\left[ -v_{i} +I( w_{j}) +\sigma \left(\overline{N}^{\alpha }\right)\right]^{2}}{2a}},\quad M_{i+1/2,j}^{\overline{N}^{\alpha }} =\frac{2M_{i,j}^{\overline{N}^{\alpha }} M_{i+1,j}^{\overline{N}^{\alpha }}}{M_{i,j}^{\overline{N}^{\alpha }} +M_{i+1,j}^{\overline{N}^{\alpha }}}.
\end{equation}

Now to complete the scheme, the only thing that remains is to specify the expression of $\overline N^\alpha$. As we shall see, the choice of $\overline N^\alpha$ will highly affect the implementations and properties of the scheme. Here we propose two choices of $\overline N^\alpha$: either $ \overline{N}^{m}$ or $ \overline{N}^{m+1}$. We call the scheme (\ref{a9OxraVzxdb})\textasciitilde (\ref{aqSw71E7wnb}) with $$ \overline{N}^{\alpha } =\overline{N}^{m}$$ the $ v$-wise semi-implicit ($ v$-SI) scheme, and the scheme (\ref{a9OxraVzxdb})\textasciitilde (\ref{aqSw71E7wnb}) with $$ \overline{N}^{\alpha } =\overline{N}^{m+1}$$ the $ v$-wise fully-implicit ($ v$-FI) scheme. The advantage of the $v$-SI scheme is that its implementation is free of nonlinear solvers, and thus is the main focus of this work. The $v$-FI scheme is mainly introduced for comparison.

To avoid ambiguity, when discussing the  $v$-SI scheme, we add a superscript "SI" for all the variables, which will be put at the first place and separated with other superscripts by a semicolon:
\begin{equation*}
    p^{m}_{i, j}, N^m_j, \overline{N}^m, H^m_j \quad\text{are written as}\quad p^{\text{SI}; m}_{i, j}, N^{\text{SI}; m}_j, \overline{N}^{\text{SI}; m}, H^{\text{SI}; m}_j \quad\text{ for the }v\text{-SI scheme};
\end{equation*}
and a similar set of notations will be introduced for the $v$-FI scheme:
\begin{equation*}
    p^{m}_{i, j}, N^m_j, \overline{N}^m, H^m_j \quad\text{are written as}\quad p^{\text{FI}; m}_{i, j}, N^{\text{FI}; m}_j, \overline{N}^{\text{FI}; m}, H^{\text{FI}; m}_j \quad \text{for the }v\text{-FI scheme}.
\end{equation*}
We still use $p^{m}_{i, j}, N^m_j, \overline{N}^m, H^m_j$ when discussing the properties that the $v$-SI scheme and the $v$-FI scheme have in common (like the conservation and positivity preserving property).

It is worth noting that, in (\ref{aiChucap+Ob}), the flux-shift part is $\overline N ^{m+1}$ for both the $v$-SI and $v$-FI construction of the scheme. This is an apparently small but important difference from the scheme proposed in \cite{JH19}. This difference will be proved to make the scheme (unconditionally) positivity preserving for all $\varepsilon$ and $\Delta v$.

\subsection{The matrix problem involved in the schemes and the iterative methods for the $ v$-FI scheme}\label{subsec: matrix form}

In this subsection, we introduce the matrix equation one will face when implementing the proposed schemes. To begin with, we define the class of matrices which will frequently appear throughout our work:
\begin{definition} \label{MatrixM}
For a given $ \overline{N} \in \mathbb{R}$ and $j=0,1,\cdots,n_w$, we define the matrix $ \boldsymbol{M}_{j}^{\overline{N}}$ as an $ n_v\times n_v$ matrix such that
\begin{equation} \label{a5A86LEWYt}
\left(\boldsymbol{M}_{j}^{\overline{N}}\right)_{kl} =\begin{cases}
-M_{k-3/2,j}^{\overline{N}} /M_{k-2,j}^{\overline{N}} & l=k-1\\
-M_{k-1/2,j}^{\overline{N}} /M_{k,j}^{\overline{N}} & l=k+1\\
\left( M_{k-3/2,j}^{\overline{N}} +M_{k-1/2,j}^{\overline{N}}\right) /M_{k-1,j}^{\overline{N}} & l=k\ \text{and} \ l\neq 1,n_v\\
M_{n-3/2,j}^{\overline{N}} /M_{n-1,j}^{\overline{N}} +1 & l=k=n_v\\
M_{1/2,j}^{\overline{N}} /M_{0,j}^{\overline{N}} & l=k=1\\
-1 & k=r+1\ \text{and} \ l=n_v\\
0 & \text{otherwise}
\end{cases}
\end{equation}
where
\begin{equation} \label{ayIywC5SuDb}
M_{i,j}^{\overline{N}} =\mathrm{e}^{-\frac{[ -v_{i} +I( w_{j}) +\sigma (\overline{N})]^{2}}{2a}}\text{ \ and \ } M_{i+1/2,j}^{\overline{N}} =\frac{2M_{i,j}^{\overline{N}} M_{i+1,j}^{\overline{N}}}{M_{i,j}^{\overline{N}} +M_{i+1,j}^{\overline{N}}} .
\end{equation}
\end{definition}
$ \boldsymbol{M}_{j}^{\overline{N}}$ is tridiagonal except for one entry at $ k=r+1$, $ l=n_v$, and its properties will be further discussed in section \ref{subsec: preparations}. For now, having defined $ \boldsymbol{M}_{j}^{\overline{N}}$, we can write the matrix form of the proposed schemes. 

For the $v$-SI scheme, assume that we have already solved $ p_{i,j}^{\text{SI} ;*m}$ from the explicit convection step (\ref{aCWU7qDHkqb}), then one easily checks\footnote{The only notable detail here is that $ N_{j}^{m+1}$ in (\ref{aiChucap+Ob}) should be replaced with $ ap_{n_{v} -1,j}^{m+1} /\Delta v$ according to (\ref{a5KFaVLjKZb}).} that solving the $v$-wise transport step (\ref{ayUOikHTLnb}) means solving the system
\begin{equation} \label{aiNYQUmvyVb}
\left( \varepsilon \boldsymbol{I} +\frac{a\Delta t}{( \Delta v)^{2}}\boldsymbol{M}_{j}^{\overline{N}^{\text{SI} ;m}}\right)\boldsymbol{p}_{:,j}^{\text{SI} ;m+1} =\varepsilon \boldsymbol{p}_{:,j}^{\text{SI} ;*m}
\end{equation}
for each $ j=0,1,\cdots ,n_{w}$, where
\begin{equation*}
\begin{aligned}
 \begin{array}{l}
\boldsymbol{p}_{:,j}^{\text{SI} ;m+1}\\
\end{array} & =\left( p_{0,j}^{\text{SI} ;m+1} ,p_{1,j}^{\text{SI} ;m+1} ,\cdots ,p_{n_{v} -1,j}^{\text{SI} ;m+1}\right)^{\top } ,\\
\boldsymbol{p}_{:,j}^{\text{SI} ;*m} & =\left( p_{0,j}^{\text{SI} ;*m} ,p_{1,j}^{\text{SI} ;*m} ,\cdots ,p_{n_{v} -1,j}^{\text{SI} ;*m}\right)^{\top } ,
\end{aligned}
\end{equation*}
$ \boldsymbol{I}$ is the $ n\times n$ identity matrix and $ \boldsymbol{M}_{j}^{\overline{N}^{\text{SI} ;m}}$ follows the definition given in (\ref{a5A86LEWYt}). (\ref{aiNYQUmvyVb}) is a linear equation for $\boldsymbol{p}_{:,j}^{\text{SI} ;m+1}$.

For the $v$-FI scheme, solving the $v$-wise transport equation means solving the system
\begin{equation} \label{au0ry5BWiV}
\left( \varepsilon \boldsymbol{I} +\frac{a\Delta t}{( \Delta v)^{2}}\boldsymbol{M}_{j}^{\overline{N}^{\text{FI} ;m+1}}\right)\boldsymbol{p}_{:,j}^{\text{FI} ;m+1} =\varepsilon \boldsymbol{p}_{:,j}^{\text{FI} ;*m}
\end{equation}
for each $ j=0,1,\cdots ,n_{w}$, where
\begin{equation*}
\begin{aligned}
 \begin{array}{l}
\boldsymbol{p}_{:,j}^{\text{FI} ;m+1}\\
\end{array} & =\left( p_{0,j}^{\text{FI} ;m+1} ,p_{1,j}^{\text{FI} ;m+1} ,\cdots ,p_{n_{v} -1,j}^{\text{FI} ;m+1}\right)^{\top } ,\\
\boldsymbol{p}_{:,j}^{\text{FI} ;*m} & =\left( p_{0,j}^{\text{FI} ;*m} ,p_{1,j}^{\text{FI} ;*m} ,\cdots ,p_{n_{v} -1,j}^{\text{FI} ;*m}\right)^{\top }.
\end{aligned}
\end{equation*}
(\ref{au0ry5BWiV}) is a  nonlinear system since the matrix $ \boldsymbol{M}_{j}^{\overline{N}^{\text{FI} ;m+1}}$ depends directly on the unsolved variable $ \overline N^{\text{FI} ;m+1}$. Therefore, we have to solve such a system with a nonlinear solver. 

We propose a fix point iteration method as follows. We set the initial guess as $ p_{i,j}^{( 0)} =p_{i,j}^{\text{FI} ;*m}$, hence $ \overline{N}^{( 0)} =\overline{N}^{m}$, then for $ k=0,1,2,\cdots $, let 
\begin{equation} \label{ayxO0an3qeb}
\begin{cases}
\displaystyle \left( \varepsilon \boldsymbol{I} +\frac{a\Delta t}{( \Delta v)^{2}}\boldsymbol{M}_{j}^{\overline{N}^{( k)}}\right)\boldsymbol{p}_{:,j}^{( k+1)} =\varepsilon \boldsymbol{p}_{:,j}^{\text{FI} ;m} \ \text{for} \ j=0,1,\cdots ,n_{w} & \\
\displaystyle \overline{N}^{( k+1)} =\Delta w\sum _{j=0}^{n_{v}} N^{( k+1)} =\Delta w\sum _{j=0}^{n_{v}}\frac{p_{n_{v} -1,j}^{( k+1)}}{\Delta v}, & 
\end{cases}
\end{equation}
where
\begin{equation*}
\boldsymbol{p}_{:,j}^{( k+1)} =\left( p_{0,j}^{( k+1)} ,p_{1,j}^{( k+1)} ,\cdots ,p_{n_{v} -1,j}^{( k+1)}\right)^{\top } .
\end{equation*}
We iterate $K$ rounds to satisfy the stopping criterion $ |\overline{N}^{(K+1)} -\overline{N}^{(K)} | \leq \epsilon$ for some small enough $\epsilon$, and then take
\begin{equation*}
\boldsymbol{p}_{:,j}^{\text{FI} ;m+1} =\boldsymbol{p}_{:,j}^{( K)} ,\ N_{j}^{\text{FI} ;m+1} =N_{j}^{( K)} ,\ \overline{N}^{\text{FI} ;m+1} =\overline{N}^{( K)}.
\end{equation*}

\section{The Properties of the Scheme} \label{sec:properties}
Having proposed the scheme in Section \ref{subsec:scheme}, we discuss how it preserves the properties of the original model.

We show that both the $v$-SI and the $v$-FI schemes proposed in Section \ref{subsec:scheme}  are mass-conservative and positivity-preserving as long as the $w$-wise CFL condition is satisfied (no matter how small $\varepsilon$ and $\Delta v$ is). 
 As for the AP property, the $v$-FI scheme and the $v$-SI scheme preserve the asymptotic limit in slightly different forms.

\subsection{Conservative properties}\label{subsec:conservation}

Being total mass invariant is one of the most basic properties of our model (\ref{auGi4Tx7zhb}), and our scheme satisfies the conservation property naturally thanks to the finite volume construction.
\begin{thm}
The schemes (\ref{a9OxraVzxdb})\textasciitilde (\ref{aqSw71E7wnb}) satisfies %\jh{I suggest stating the conservation without using $p^*$ (you may use it in the proof).}
\begin{equation} \label{aC6l+BodKtb}
\sum _{i=0}^{n_{v}}\sum _{j=0}^{n_{w}} p_{i,j}^{m+1} =\sum _{i=0}^{n_{v}}\sum _{j=0}^{n_{w}} p_{i,j}^{m},
\end{equation}
\end{thm}
\blue{Proof: Summing (\ref{a9OxraVzxdb}) over $ i$ from $ 0$ to $ n_{v} -1$, over $ j$ from $ 0$ to $ n_{w}$ and using (\ref{aOpH3A49bX}) (\ref{aiChucap+Ob}), we have}
\blue{\begin{equation}\label{34irgodg}
\begin{aligned}
&\sum_{i=0}^{n_v-1}\sum_{j=0}^{n_w} p_{i,j}^{m+1} \\= & \sum_{i=0}^{n_v-1}\sum_{j=0}^{n_w} p_{i,j}^{m} + \frac{\Delta t}{\varepsilon \Delta v} \sum_{i=0}^{n_v-1}\sum_{j=0}^{n_w} (F^m_{i+1/2,j}-F^m_{i-1/2,j}) -  \frac{\Delta t}{\Delta w}\sum_{i=0}^{n_v-1}\sum_{j=0}^{n_w} (\Psi^m_{i,j+1/2} - \Psi^m_{i,j-1/2}) \\
= & \sum_{i=0}^{n_v-1}\sum_{j=0}^{n_w} p_{i,j}^{m} + \frac{\Delta t}{\varepsilon \Delta v} \sum_{j=0}^{n_w} (F^m_{n_v-1/2,j}-F^m_{-1/2,j}) - \frac{\Delta t}{\Delta w} \sum_{i=0}^{n_v-1} (\Psi^m_{i,n_w+1/2} - \Psi^m_{i,-1/2}) \\
= & \sum_{i=0}^{n_v-1}\sum_{j=0}^{n_w} p_{i,j}^{m} + \frac{\Delta t}{\varepsilon \Delta v} \sum_{j=0}^{n_w} (0-0) - \frac{\Delta t}{\Delta w} \sum_{i=0}^{n_v-1} (0 - 0) = \sum_{i=0}^{n_v-1}\sum_{j=0}^{n_w} p_{i,j}^{m},
\end{aligned}
\end{equation}}
\blue{and because of the boundary condition $ p_{n_{v} ,j}^{m} =0$, we can change ``$\sum_{i=0}^{n_v-1}$'' on both ends of (\ref{34irgodg}) into ``$\sum_{i=0}^{n_v}$'', which gives us (\ref{aC6l+BodKtb}).}

\subsection{Properties of the matrix \(M_{j}^{\overline N}\)}
\label{subsec: preparations}

Before discussing the properties of our schemes, we firstly investigate the matrix $ \boldsymbol{M}_{j}^{\overline{N}}$ defined in (\ref{a5A86LEWYt}) and (\ref{ayIywC5SuDb}). \blue{This matrix has a one dimensional kernel spanned by a strictly positive vector}. This property will help us to show the positivity of the schemes.

The main result in this section is
\begin{lema} \label{lema: positive 1D kernel}
For all $ \overline{N} \geq 0$, $ \boldsymbol{M}_{j}^{\overline{N}}$ defined in (\ref{a5A86LEWYt}) and (\ref{ayIywC5SuDb})  has a unique (up to a positive multiple) right eigenvector $ \boldsymbol{q} =( q_{0} ,\cdots ,q_{n_{v} -1})^{\top }$ for eigenvalue $ 0$ which is strictly positive. %(except for positive multiples) \jh{Do you mean unique up to a positive multiple?}.
\end{lema}
Proof: We want to study the solution of
\begin{equation} \label{aCa2ap+sZL}
{\boldsymbol M_{j}^{\overline{N}} \boldsymbol q} = \boldsymbol 0.
\end{equation}
We left-multiply both sides of (\ref{aCa2ap+sZL}) by the $ n\times n$ bi-diagonal matrix
\begin{equation*}
\boldsymbol{L} =\begin{bmatrix}
1 &  &  &  & \\
1 & 1 &  &  & \\
 & 1 & 1 &  & \\
 &  & \ddots  & \ddots  & \\
 &  &  & 1 & 1
\end{bmatrix} ,\text{ or entry-by-entry} :\ L_{kl} =\begin{cases}
1 & k=l+1\ \text{or} \ k=l\\
0 & \text{otherwise}
\end{cases}
\end{equation*}
obtaining $ \left(\boldsymbol{LM}_{j}^{\overline{N}}\right)\boldsymbol{q} =\boldsymbol 0$ (note that $\boldsymbol L$ is non-singular, so $\boldsymbol L \boldsymbol M_j^{\overline N}$ has the same kernel with $\boldsymbol M_j^{\overline N}$), which can be rearranged into the following system: %\jh{I don't quite understand this notation. What does (3.4) represent?} \hq{(3.4) is the equation system version of $ \left(\boldsymbol{LM}_{j}^{\overline{N}}\right)\boldsymbol{q} =\boldsymbol 0$, from which we can derive that all $q_i$'s $(i=0,\cdots,n_v-1)$ have the same sign}
%\hq{``$q_i=$" is moved to the right of the brackets}
\begin{equation} \label{ayACrorIpkb}
\begin{cases}
\displaystyle q_{i} =\frac{M_{i,j}^{\overline{N}}}{M_{i,j+1}^{\overline{N}}} q_{i+1} & i=r-1,r-2,\cdots ,0\\[15pt]
\displaystyle q_{i} =\frac{M_{i,j}^{\overline{N}}}{M_{i,j+1}^{\overline{N}}} q_{i+1} +\frac{M_{i,j}^{\overline{N}}}{M_{i,j+1/2}^{\overline{N}}} q_{n-1} & i=n-3,n-4,\cdots ,r\\[15pt]
\displaystyle q_{i} =\left(\frac{M_{i,n-2}^{\overline{N}}}{M_{i,n-1}^{\overline{N}}} +\frac{M_{i,n-2}^{\overline{N}}}{M_{i,n-3/2}^{\overline{N}}}\right) q_{n-1} & i=n-2
\end{cases}
\end{equation}
One can derive from (\ref{ayACrorIpkb}) recurrently (in a reversed order, from $ i=n-2$ to $ i=0$) that $ q_{n-2} ,q_{n-3} ,\cdots, q_{0}$ can all be written as a positive multiple of $ q_{n-1}$ whose expressions only contains entries of $ \boldsymbol{M}_{j}^{\overline{N}}$. This implies strict-positiveness and uniqueness (up to a positive number multiplication) of $ \boldsymbol{q}$.\hfill$\square$

\subsection{Positivity preserving properties}\label{subsec: positivity preserving}

Solutions being non-negative is another basic property of model (\ref{auGi4Tx7zhb}), hence we also expect our scheme to be positivity preserving. In Theorem 2.2 \cite{JH19}, it is proved that the scheme's positivity preserving property relies on the grid ratio $ \Delta t/( \Delta v)^{2}$. In this section, we see that the positivity preserving property of our scheme does not rely on $ \Delta t/( \Delta v)^{2}$ or $\varepsilon$, and this improvement comes from our implicit treatment to the flux-shift operator that is stated in Section \ref{sec: The schemes}. %\jh{Section 2? Also, it helps to mention the parabolic CFL condition in the discussion before section 2.3.}

We start with the general matrix form of our problem (both for the $ v$-SI scheme and for the $ v$-FI scheme):
\begin{equation} \label{aam-SykHDY}
\left( \varepsilon \boldsymbol{I} +\frac{a\Delta t}{( \Delta v)^{2}}\boldsymbol{M}_{j}^{\overline{N}^{\alpha }}\right)\boldsymbol{p}_{:,j}^{m+1} =\varepsilon \boldsymbol{p}_{:,j}^{*m},
\end{equation}
where
\begin{equation*}
\begin{aligned}
 \begin{array}{l}
\boldsymbol{p}_{:,j}^{m+1}\\
\end{array} & =\left( p_{0,j}^{m+1} ,p_{1,j}^{m+1} ,\cdots ,p_{n_{v} -1,j}^{m+1}\right)^{\top } ,\\
\boldsymbol{p}_{:,j}^{*m} & =\left( p_{0,j}^{*m} ,p_{1,j}^{*m} ,\cdots ,p_{n_{v} -1,j}^{*m}\right)^{\top } .
\end{aligned}
\end{equation*}
A crucial property of system (\ref{aam-SykHDY}) is that the to-be-inverted matrix $ \varepsilon \boldsymbol{I} +\frac{a\Delta t}{( \Delta v)^{2}}\boldsymbol{M}_{j}^{\overline{N}^{\alpha }}$ is the so-called M-matrix, which is an important topic in the area of positive operators. In \cite{ple77}, the authors provided forty equivalent statements that ensures a Z-matrix (a matrix with all diagonal entries non-negative and all non-diagonal entries non-positive) to be an M-matrix. We pick out the statement F16 and K33 in this article, forming the following lemma: %\jh{Bad sentence here.}
\begin{lema}[\red{Equivalence of statements F16 and K33 of Theorem 1 in \cite{ple77}}]
A semi-positive Z-matrix is an M-matrix, and henceforth monotone. Specifically, if a $ n\times n$ matrix $ \boldsymbol{M}$ satisfies
\begin{equation*}
\begin{cases}
\boldsymbol{M}_{i,j} \geq 0 & j=i\\
\boldsymbol{M}_{i,j} \leq 0 & j\neq i
\end{cases}
\end{equation*}
and there exists a strictly positive (every entry is positive) $ n\times 1$ vector $ \boldsymbol{v}^{*}  >0$ such that $ \boldsymbol{Mv}^{*}  >0$, %then for all non-negative \ $ n\times 1$ vector $ v$ \jh{what do you mean here?}, 
then $ \boldsymbol{M}$ is a non-singular M-matrix, and henceforth $ \boldsymbol{Mv} \geq 0$ implies $ \boldsymbol{v} \geq 0$. Here 

- $ \boldsymbol{v}  >0$ means that all the entries of $ \boldsymbol{v}$ are positive, and

- $ \boldsymbol{v} \geq 0$ means that all the entries of $ \boldsymbol{v}$ are non-negative.
\label{lema:ple77}
\end{lema}
With Lemma \ref{lema:ple77}, we can further prove the following Lemma:
\begin{lema}
For (\ref{aam-SykHDY}), $ \boldsymbol{p}_{:, j}^{*m} \geq 0$ indicates $ \boldsymbol{p}_{:, j}^{m+1} \geq 0$.
\end{lema}
\textit{Proof}: Firstly, it is not hard to check that $ \varepsilon \boldsymbol{I} +\frac{a\Delta t}{( \Delta v)^{2}}\boldsymbol{M}_{j}^{\overline{N}^{\alpha}}$ is a Z-matrix according to the definitions. Secondly, we notice that
\begin{equation*}
\left( \varepsilon \boldsymbol{I} +\frac{a\Delta t}{( \Delta v)^{2}}\boldsymbol{M}_{j}^{\overline{N}^{\alpha }}\right)\boldsymbol{Q}_{:,j}^{\overline{N}^{\alpha }} =\varepsilon {\boldsymbol Q_{:,j}^{\overline{N}^{\alpha }}}  >0
\end{equation*}
since ${\boldsymbol Q_{i,j}^{\overline{N}^{\alpha}}}$ is in the kernel of $ \boldsymbol{M}_{j}^{\overline{N}^{\alpha}}$, according to Lemma 2 (taking $ \boldsymbol{v}^{*}$ as $ \boldsymbol Q_{j}^{\overline{N}^{\alpha}}$), the matrix $ \varepsilon \boldsymbol{I} +\frac{a\Delta t}{( \Delta v)^{2}}\boldsymbol{M}_{j}^{\overline{N}^{\alpha}}$ is a non-singular M-Matrix and henceforth \ $ \boldsymbol{p}_{j}^{*m} \geq 0$ indicates $ \boldsymbol{p}_{j}^{m+1} \geq 0$.\hfill$\square$

Applying Lemma 3 and (\ref{aCWU7qDHkqb}), we can give out the following sufficient conditions of the positivity preserving property of the proposed schemes
\begin{thm}
(Positivity Preserving Properties) \ $ p_{i,j}^{m+1} \geq 0,$ holds for all $ i=0,\cdots ,n_{v}$ and $ j=0,\cdots ,n_{w}$ as long as
\begin{equation} \label{amSYbT13I+}
p_{i,j}^{m} -\frac{\Delta t}{\Delta w}\left( \Phi _{i,j+1/2}^{m} -\Phi _{i,j-1/2}^{m}\right) \geq 0
\end{equation}
holds for all $ i=0,\cdots ,n_{v} ,j=0,\cdots ,n_{w}$.
\end{thm}

(\ref{amSYbT13I+}) tells us that as long as the grid ratio $ \Delta t/\Delta w$ is not too big, (which is the most basic requirement for almost all numerical schemes for hyperbolic conservation law, and is also known as the CFL condition), the numerical solution will be non-negative. Such a property is unrelated to the grid ratio $\Delta t / \Delta v^2$ or the scaling parameter $\varepsilon$. This means that we can take arbitrarily small $\Delta v$ and $\varepsilon$ without worrying about violations of the positivity of the solutions.

Especially, as we can take $\varepsilon$ arbitrarily small for a fixed $\Delta t$, it becomes possible for us to futher analyze the asymptotic behavior of the scheme when $\varepsilon\to 0$.

\subsection{Asymptotic preserving properties} \label{subsec:asymptotic preserving}

% \jh{We don't have to repeat equations. Rather it is better to mention that we first discuss the asymptotic properties of the original/continuous model.}
We first discuss the asymptotic properties of the original continuous model. (\ref{auGi4Tx7zhb})%. We rewrite equations (\ref{auGi4Tx7zhb}) and (\ref{aeRiM-jHclb}) here:
% \begin{equation} \label{aGVz7d-WdM}
% \begin{cases}
% \begin{aligned}
% \frac{\partial p}{\partial t} +\frac{\partial }{\partial w} [(\overline{N} (t)N(w,t)K(w)-w)p] & \\
% =\frac{1}{\varepsilon }\Bigl\{a\frac{\partial ^{2} p}{\partial v^{2}} -\frac{\partial }{\partial v} [(-v+I(w) & +w\sigma (\overline{N} (t)))p]+N(w,t)\delta (v-V_{R} )\Bigr\}\\
% \text{for } (v,w,t) & \in (-\infty ,V_{F} )\times (-\infty ,+\infty )\times (0,T)
% \end{aligned}\\
% {\displaystyle N(w,t):=-a\frac{\partial p}{\partial v} (V_{F} -0,w,t) >0,\ \ \overline{N} (t)=\int _{-\infty }^{\infty } N(w,t)\mathrm{d} w.}
% \end{cases}
% \end{equation}
% and
% \begin{equation} \label{aGBOLGAkSPb}
% \int _{-\infty }^{V_{F}} p( v,w,t)\mathrm{d} v=H( w,t).
% \end{equation}

Integrating the first equation of (\ref{auGi4Tx7zhb}) over $\displaystyle v$ from $\displaystyle -\infty $ to $\displaystyle V_{F}$ and applying the boundary condition, we get
\begin{equation} \label{aGzCbxmES9}
\frac{\partial H}{\partial t} +\frac{\partial }{\partial w}[ (\overline{N} (t)N(w,t)K(w)-w)H] =0,
\end{equation}
where $H=H(w,t)$ is defined in (\ref{aeRiM-jHclb}). When taking $\displaystyle \varepsilon \rightarrow 0$ in (\ref{auGi4Tx7zhb}), we have
\begin{equation} \label{ai8vhCSh+U}
\begin{cases}
\displaystyle a\frac{\partial ^{2} p}{\partial v^{2}} -\frac{\partial }{\partial v} [(-v+I(w)+w\sigma (\overline{N} (t)))p]+N(w,t)\delta (v-V_{R} )=0, & \\
{\displaystyle N(w,t):=-a\frac{\partial p}{\partial v} (V_{F} -0,w,t),\ \ \overline{N} (t)=\int _{-\infty }^{\infty } N(w,t)\mathrm{d} w.} & 
\end{cases}
\end{equation}
%\hq{I clarified here that the steady state solution exists.}
According to Theorem 3.1 in \cite{BP18}, the solution of equation (\ref{ai8vhCSh+U}) is graranteed to exist and be unique when $\displaystyle H( w,t)$ is given.
Equations (\ref{aeRiM-jHclb}), (\ref{aGzCbxmES9}) and (\ref{ai8vhCSh+U}) gives the asymptotic limit equation of (\ref{auGi4Tx7zhb}) as $\displaystyle \varepsilon \rightarrow 0$, where (\ref{aGzCbxmES9}) can be viewed as the macroscopic equation governing the slow dynamics while (\ref{ai8vhCSh+U}) gives the local microscopic equilibrium of the fast dynamics. %\jh{I recall there is a complicated fixed point argument for the existence of the asymptotic state. Do we want to mention that?} \hq{I gave up mentioning that. I think it will only make the simple problem complicated if we describe the scheme in a semi-discretized form}
We consider our scheme with only time $\displaystyle t$ discretized
\begin{equation} \label{a9ucjp0eoQb}
\begin{cases}
\begin{aligned}
\frac{p^{n+1} -p^{n}}{\Delta t} +\frac{\partial }{\partial w} [(\overline{N}^{n} N^{n} (w)K(w)-w) & p^{n} ]\\
=\frac{1}{\varepsilon }\Bigl\{a\frac{\partial ^{2} p^{n+1}}{\partial v^{2}} -\frac{\partial }{\partial v} [(-v+I(w)+ & w\sigma (\overline{N}^{\alpha } ))p^{n+1} ]+N^{n+1} (w)\delta (v-V_{R} )\Bigr\}\\
\text{for } (v,w,t) & \in (-\infty ,V_{F} )\times (-\infty ,+\infty )\times (0,T)
\end{aligned}\\
{\displaystyle N^{n+1} (w):=-a\frac{\partial p}{\partial v} (V_{F} -0,w,t),\ \ \overline{N}^{n+1} =\int _{-\infty }^{\infty } N^{n+1} (w)\mathrm{d} w;}
\end{cases}
\end{equation}
where $\displaystyle \alpha $ is $\displaystyle n$ and $\displaystyle n+1$ for the $\displaystyle v$-SI scheme and $\displaystyle v$-FI scheme, respectively.

For the $\displaystyle v$-FI scheme, integrating the first equation in (\ref{a9ucjp0eoQb}) and applying the boundary conditions, we have
\begin{equation} \label{aGmsynVXERb}
\frac{H^{n+1} -H^{n}}{\Delta t} +\frac{\partial }{\partial w}\left[ (\overline{N}^{n} N^{n} (w)K(w)-w)H^{n}\right] =0
\end{equation}
which is a discretization of (\ref{aGzCbxmES9}); and when $\displaystyle \varepsilon \rightarrow 0$, (\ref{a9ucjp0eoQb}) becomes
\begin{equation}\label{whueirvbi}
\begin{cases}
\displaystyle a\frac{\partial ^{2} p^{n+1}}{\partial v^{2}} -\frac{\partial }{\partial v} [(-v+I(w)+w\sigma (\overline{N}^{n+1} ))p^{n+1} ]+N^{n+1} (w)\delta (v-V_{R} )=0 & \\
{\displaystyle N^{n+1} (w):=-a\frac{\partial p}{\partial v} (V_{F} -0,w,t),\ \ \overline{N}^{n+1} =\int _{-\infty }^{\infty } N^{n+1} (w)\mathrm{d} w,} & 
\end{cases}
\end{equation}
which is directly the same as (\ref{ai8vhCSh+U}).

For the $\displaystyle v$-SI scheme, Integrating the first equation in (\ref{a9ucjp0eoQb}) and applying the boundary conditions, we obtain (\ref{aGmsynVXERb}) again, while taking $\displaystyle \varepsilon \rightarrow 0$ in (\ref{a9ucjp0eoQb}), we have
\begin{equation} \label{aCTcNR2-Sl}
\begin{cases}
\displaystyle a\frac{\partial ^{2} p^{n+1}}{\partial v^{2}} -\frac{\partial }{\partial v} [(-v+I(w)+w\sigma (\overline{N}^{n} ))p^{n+1} ]+N^{n+1} (w)\delta (v-V_{R} )=0 & \\
{\displaystyle N^{n+1} (w):=-a\frac{\partial p}{\partial v} (V_{F} -0,w,t),\ \ \overline{N}^{n+1} =\int _{-\infty }^{\infty } N^{n+1} (w)\mathrm{d} w.} & 
\end{cases}
\end{equation}
We remark that the local equilibrium (\ref{aCTcNR2-Sl}) can be viewed as a linearization of (\ref{whueirvbi}) such that the nonlinear coefficient $\sigma(\overline N^{n+1})$ is replaced by $\sigma(\overline N^{n})$,

The schemes (\ref{aGmsynVXERb}) and (\ref{aCTcNR2-Sl}) can be viewed as the limiting scheme for the following time rescaled maybe ``delayed" equation as $\varepsilon \rightarrow 0$
\begin{equation} \label{aqKc9bZluSb}
\begin{cases}
\begin{aligned}
\frac{\partial p}{\partial t} +\frac{\partial }{\partial w} [(\overline{N} (t)N(w,t)K(w) & -w)p]\\
=\frac{1}{\varepsilon }\Bigl\{a\frac{\partial ^{2} p}{\partial v^{2}} -\frac{\partial }{\partial v} [(-v+ & I(w)+w\sigma (\overline{N} (t-\Delta t)))p]+N(w,t)\delta (v-V_{R} )\Bigr\}\\
\text{for } (v,w,t) & \in (-\infty ,V_{F} )\times (-\infty ,+\infty )\times (0,T)
\end{aligned}\\
{\displaystyle N (w):=-a\frac{\partial p}{\partial v} (V_{F} -0,w,t),\ \ \overline{N} =\int _{-\infty }^{\infty } N (w)\mathrm{d} w.}
\end{cases}
\end{equation}
Therefore, the $\displaystyle v$-SI scheme is asymptotic preserving up to a $\displaystyle \Delta t$ time delay. In other words, when $\varepsilon\rightarrow 0$, it is still a consistent first-order in time discretization to the limit of the original model.

\section{Numerical Results}
\label{sec:numerical_results}
This section gives the numerical results obtained from our schemes, involving both the tests on the schemes' performance and the explorations of the model by our scheme. In Section \ref{subsec:convergence}, we numerically test the order of accuracy of our schemes; in Section \ref{subsec: numtest AP}, we numerically validate the asymptotic preserving properties of our schemes; in Section \ref{subsec: Learning and Reacting}, we test the model’s learning and discrimination abilitites; in Section \ref{subsec: Positive Synaptic Weights}, we focus on exploring the model’s behavior when the neural system is excitatory.
\subsection{Order of Convergence}\label{subsec:convergence}
In this part, we test the order of accuracy of the $v$-SI schemes. Since the exact solution is unavailable, we estimate the order of the error by
$$
O_{h, L^p}=\log _{2} \frac{\left\|\omega_{h}-\omega_{\frac{h}{2}}\right\|_p}{\left\|\omega_{\frac{h}{2}}-\omega_{\frac{h}{4}}\right\|_p},
$$
where $\omega_{h}$ is the numerical solution with step length $h$. The term $O_{h}$ above is an approximation for the accuracy order. Errors in both $L^1$ and $L^2$ norms are examined.

We choose $V_F=2, V_R=1, V_{\min} = -4, a = 1, \varepsilon = 0.5, W_{\min} = -1.1, W_{\max} = 0.1, \sigma(\overline N) = \overline N$, $I(w) \equiv 0$ and
$$
p_{\text{init}} = \begin{cases}
\displaystyle \sin ^2(\pi v) \sin ^2 (\pi w) & -1 < w < 0 \text{ and } -1 < v < 1\\
0 & \text{otherwise}
\end{cases}
$$
throughout this section, and we test the accuracy for both $T_{\max } = 0.1$ and $T_{\max } = 1$. We test the numerical results on $v, w, t$ directions by fixing two of $\Delta t, \Delta w$ and $\Delta v$ while adjusting the third one. The results are shown in Table \ref{tab: v-ooa}, Table \ref{tab: w-ooa} and Table \ref{tab: t-ooa}, respectively.

\begin{table}[!ht]
\begin{center}

\begin{tabular}{|c|c|c|c|c|}
\hline 
 $\displaystyle \Delta v$ & $\displaystyle \| p_{\Delta v} -p_{\Delta v/2} \| _{1}$ & $\displaystyle O_{\Delta v,L^{1}}$ & $\displaystyle \| p_{\Delta v} -p_{\Delta v/2} \| _{2}$ & $\displaystyle O_{\Delta v,L^{2}}$ \\
\hline 
 2.0e-1 & 1.1337e-03 & 2.0818 & 3.5479e-03 & 2.0675 \\
\hline 
 1.0e-1 & 2.6781e-04 & 2.0122 & 8.4643e-04 & 2.0080 \\
\hline 
 5.0e-2 & 6.6390e-05 & 1.9340 & 2.1044e-04 & 1.8739 \\
\hline 
 2.5e-2 & 1.7374e-05 & - & 5.7417e-05 & - \\
 \hline
\end{tabular}

\vspace{10pt}

\begin{tabular}{|c|c|c|c|c|}
\hline 
 $\displaystyle \Delta v$ & $\displaystyle \| p_{\Delta v} -p_{\Delta v/2} \| _{1}$ & $\displaystyle O_{\Delta v,L^{1}}$ & $\displaystyle \| p_{\Delta v} -p_{\Delta v/2} \| _{2}$ & $\displaystyle O_{\Delta v,L^{2}}$ \\
\hline 
 2.0e-1 & 1.4571e-01 & 1.8905 & 2.8929e-01 & 1.7984 \\
\hline 
 1.0e-1 & 3.9299e-02 & 1.7441 & 8.3168e-02 & 1.7309 \\
\hline 
 5.0e-2 & 1.1731e-02 & 1.7872 & 2.5056e-02 & 1.8471 \\
\hline 
 2.5e-2 & 3.3990e-03 & - & 6.9641e-03 & - \\
 \hline
\end{tabular}
\end{center}
\caption{$v$-wise order of accuracy, $\displaystyle \Delta t=10^{-3} ,\ \Delta w= 10^{-2}$ are fixed. Upper table: $T_{\max } = 0.1$; Lower table: $T_{\max } = 2.5$.}
\label{tab: v-ooa}
\end{table}

\begin{table}[!ht]
\begin{center}

\begin{tabular}{|c|c|c|c|c|}
\hline 
 $\displaystyle \Delta w$ & $\displaystyle \| p_{\Delta w} -p_{\Delta w/2} \| _{1}$ & $\displaystyle O_{\Delta w,L^{1}}$ & $\displaystyle \| p_{\Delta w} -p_{\Delta w/2} \| _{2}$ & $\displaystyle O_{\Delta w,L^{2}}$ \\
\hline 
 4.0e-2 & 4.2858e-03 & 0.9550 & 9.9587e-03 & 0.9543 \\
\hline 
 2.0e-2 & 2.2108e-03 & 1.0038 & 5.1394e-03 & 1.0030 \\
\hline 
 1.0e-2 & 1.1025e-03 & 0.9849 & 2.5643e-03 & 0.9801 \\
\hline 
 0.5e-2 & 5.5703e-04 & - & 1.3000e-03 & - \\
 \hline
\end{tabular}

\vspace{10pt}

\begin{tabular}{|c|c|c|c|c|}
\hline 
 $\displaystyle \Delta w$ & $\displaystyle \| p_{\Delta w} -p_{\Delta w/2} \| _{1}$ & $\displaystyle O_{\Delta w,L^{1}}$ & $\displaystyle \| p_{\Delta w} -p_{\Delta w/2} \| _{2}$ & $\displaystyle O_{\Delta w,L^{2}}$ \\
\hline 
 4.0e-2 & 4.5348e-01 & 1.0102 & 3.9149e-01 & 0.7353 \\
\hline 
 2.0e-2 & 2.2514e-01 & 1.0095 & 2.3516e-01 & 0.7557 \\
\hline 
 1.0e-2 & 1.1183e-01 & 0.9717 & 1.3928e-01 & 0.6145 \\
\hline 
 0.5e-2 & 5.7022e-02 & - & 9.0971e-02 & - \\
 \hline
\end{tabular}
\end{center}
\caption{$w$-wise order of accuracy, $\Delta t=10^{-3} ,\ \Delta v=10^{-1}$ are fixed. Upper table: $T_{\max } = 0.1$; Lower table: $T_{\max } = 2.5$.}
\label{tab: w-ooa}
\end{table}

\begin{table}[!ht]
\begin{center}

\begin{tabular}{|c|c|c|c|c|}
\hline 
 $\displaystyle \Delta t$ & $\displaystyle \| p_{\Delta t} -p_{\Delta t/2} \| _{1}$ & $\displaystyle O_{\Delta t,L^{1}}$ & $\displaystyle \| p_{\Delta t} -p_{\Delta t/2} \| _{2}$ & $\displaystyle O_{\Delta t,L^{2}}$ \\
\hline 
 2.0e-3 & 3.8523e-04 & 0.9730 & 1.2219e-03 & 0.9647 \\
\hline 
 1.0e-3 & 1.9625e-04 & 0.9686 & 6.2608e-04 & 0.9626 \\
\hline 
 5.0e-4 & 1.0028e-04 & 1.0093 & 3.2125e-04 & 1.0089 \\
\hline 
 2.5e-4 & 4.9819e-05 & - & 1.5963e-04 & - \\
 \hline
\end{tabular}

\vspace{10pt}

\begin{tabular}{|c|c|c|c|c|}
\hline 
 $\displaystyle \Delta t$ & $\displaystyle \| p_{\Delta t} -p_{\Delta t/2} \| _{1}$ & $\displaystyle O_{\Delta t,L^{1}}$ & $\displaystyle \| p_{\Delta t} -p_{\Delta t/2} \| _{2}$ & $\displaystyle O_{\Delta t,L^{2}}$ \\
\hline 
 2.0e-3 & 4.9612e-03 & 0.9473 & 4.5521e-03 & 0.9676 \\
\hline 
 1.0e-3 & 2.5729e-03 & 0.9884 & 2.3278e-03 & 0.9954 \\
\hline 
 5.0e-4 & 1.2969e-03 & 0.9781 & 1.1677e-03 & 0.9895 \\
\hline 
 2.5e-4 & 6.5836e-04 & - & 5.8810e-04 & - \\
 \hline
\end{tabular}
\end{center}
\caption{$t$-wise order of accuracy, $\Delta v=10^{-1} ,\ \Delta w=10^{-2}$ are fixed. Upper table: $T_{\max } = 0.1$; Lower table: $T_{\max } = 2.5$.}
\label{tab: t-ooa}
\end{table}

\begin{figure}
\centering   
\includegraphics[scale=0.3]{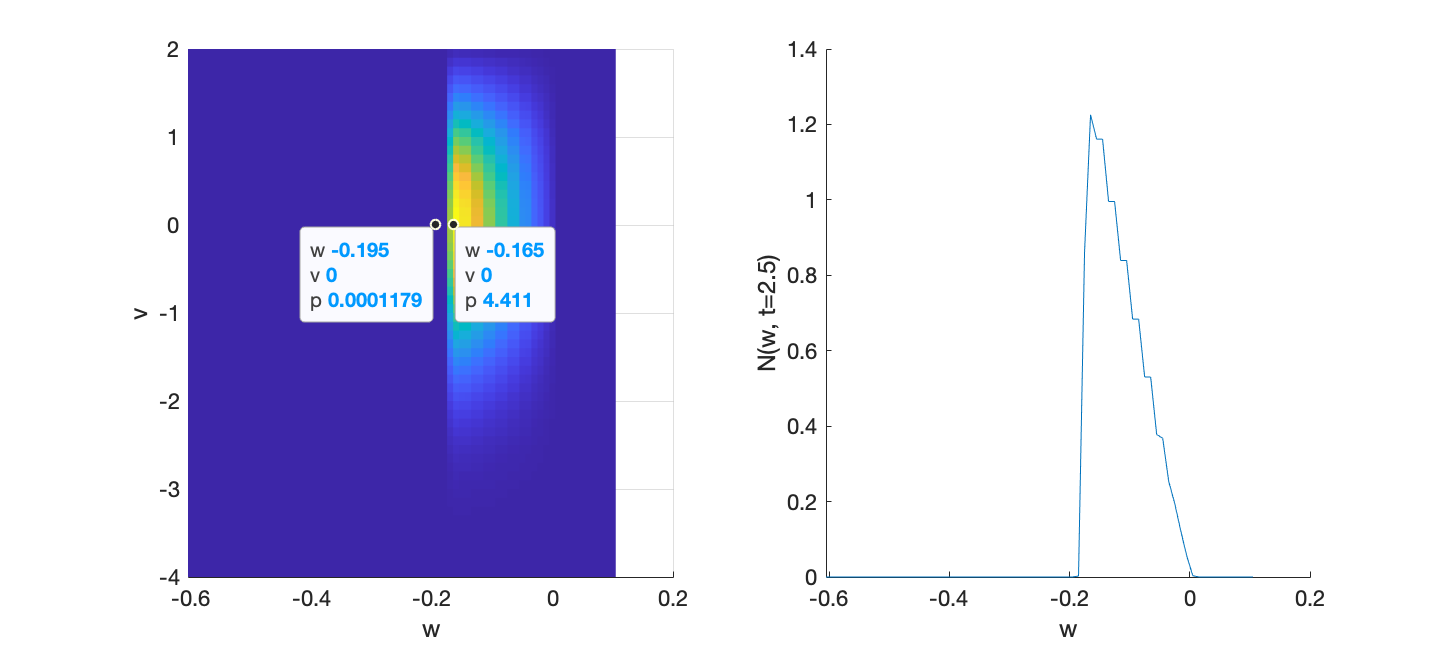}
    \caption{When $t=2.5$, there is a discontinuity in the numerical solution. Left: $p(v,w,t=2.5)$, with the marked data showing that small variance along the $w$ direction causes a big variance of $p$; Right: $N(w,t=2.5)$.}
\label{ooa_shock_show}
\end{figure}

It can be seen that when $T_{\max } = 0.1$, the scheme shows the expected  second-order convergence rate in the $v$ direction, and the first-order convergence in $w, t$ directions. However, when $T_{\max } = 2.5$, the $v$- and $w$- wise orders of accuracy decrease a bit. This may be caused by the $w$-wise nonlinearity of the equation (\ref{auGi4Tx7zhb}). Indeed, when $T=2.5$ the solution has already developed into a discontinuous pattern (Figure \ref{ooa_shock_show}).% \jh{Are you only testing the semi-implicit scheme here?}

\subsection{Asymptotic behaviors of the $v$-SI and $v$-FI schemes}
\label{subsec: numtest AP}
In this subsection, we numerically validate the asymptotic preserving properties of both the $v$-SI scheme and the $v$-FI scheme by validating that with fixed grid size, the numerical solution of the scheme converge to the $v$-wise quasi-steady state at the discretized level
%defined in Definition \ref{v-wise-QS} 
as $\varepsilon\to 0$. %\jh{Again be consistent, the asymptotics is not discussed until section 3.} \hq{This `2' is Definition number}

Firstly we introduce a discretized version of the time-independent quasi-steady state equations (\ref{aeRiM-jHclb}) and (\ref{ai8vhCSh+U}) for a given $H=H(w)$. This will serve as a reference solution. We denote the solution of the continuous time-independent quasi-steady state equations (\ref{aeRiM-jHclb}) and (\ref{ai8vhCSh+U}) for a given $H=H(w)$ as $P^H(v,w)$. Since (\ref{ai8vhCSh+U}) means that $ P^H(v,w)$ should nullify the stiff operator in (\ref{auGi4Tx7zhb}), our discretized version can also be directly defined to nullify the stiff terms (order $ 1/\varepsilon $ terms) in the main equation (\ref{a9OxraVzxdb}) in our scheme. % This gives us 
% \begin{definition} \label{v-wise-QS}
 For a given discretized grid function $H_j$ such that $\Delta w\sum _{j=0}^{n_{w}} H_{j} =1$
the $v$-wise quasi-steady state in the discretized level $ P_{i,j}^{H}$ and corresponding $ N_{j}^{H} ,\overline{N}^{H}$ is defined to satisfy the following equations %\hq{there is no longer a ``Definition" here.}
\begin{equation} \label{ay8hBjpeac}
\Delta v\sum _{i=0}^{n_{v} -1} P_{i,j}^{H} =H_{j} \ \ \ \text{for} \ j=0,1,\cdots ,n_{w} ,
\end{equation}
\begin{equation} \label{ayzeNWYrP8}
F_{i+1/2,j}^{H} -F_{i-1/2,j}^{H} =0
\end{equation}
\begin{equation} \label{aevf4jCieXb}
F_{i+1/2,j}^{H} =\begin{cases}
a\frac{M_{i+1/2,j}^{\overline{N}^{H}}}{\Delta v}\left(\frac{P_{i+1,j}^{H}}{M_{i+1,j}^{\overline{N}^{H}}} -\frac{P_{i,j}^{H}}{M_{i,j}^{\overline{N}^{H}}}\right) +N_{j}^{H}\eta ( v_{i+1/2} - V_R) & i=0,1,\cdots ,n_{v} -2\\
0 & i=-1,n_{v} -1
\end{cases} \quad 
\end{equation}
\begin{equation} \label{aKk9hL9B6y}
N_{j}^{H} =aP_{n_{v} -1,j}^{H} /\Delta t,\ \overline{N}^{H} =\Delta w\sum _{j=0}^{n_{w}} N_{j}^{H}
\end{equation}
% \end{definition}

 To generate this discretized $v$-wise quasi-steady state for a given $H_j$ with $\Delta v \sum_{j=0}^{n_w} H_j=1$, one may solve the system
\begin{equation} \label{fix-point-iteration-P}
\begin{cases}
{\displaystyle \boldsymbol{M}_{j}^{\overline{N}^{(k)}}\boldsymbol{P}_{:,j}^{(k+1)} =\boldsymbol{0} \ \text{and} \ \Delta t\sum_{i=0}^{n_v-1}{P}_{i,j}^{(k+1)} =H_{j} \ \text{for} \ j=0,1,\cdots ,n_{w}} & \\
{\displaystyle \overline{N}^{(k+1)} =\Delta w\sum _{j=0}^{n_{v}} N^{(k+1)} =\Delta w\sum _{j=0}^{n_{v}}\frac{p_{n_{v} -1,j}^{(k+1)}}{\Delta v}} & 
\end{cases}
\end{equation}
iteratively from some initial guess of ${\overline N}^{(0)}$, and take $P_{i,j}^H = P_{i,j}^{(K)}, N_{j}^H = N_{j}^{(K)}, {\overline N}^H = {\overline N}^{(K)}$ for sufficiently large $K$ that satisfies some stopping criterion like $|N^{(K+1)}-N^{(K)}|\leq \epsilon$ for some small enough $\epsilon$.

The numerical experiments are conducted as follows: 
\begin{enumerate}
\item Fix all the scheme parameters except $\varepsilon$ (these parameters include all the grid lengths, initial/boundary values and coefficients, especially, $\Delta t$ is fixed), and solve the scheme. The solution for a fixed $\varepsilon$ is denoted as $p_\varepsilon$.

\item Having solved for $p_\varepsilon$, find $H_\varepsilon$ by (\ref{ayvWHg2kk2}) %{\jh{??}
, and thereafter the corresponding discretized $v$-wise quasi-steady state $P^{H_{\varepsilon}}$ by (\ref{ay8hBjpeac})\textasciitilde (\ref{aKk9hL9B6y}), which can be numerically obtained by implementing the iterations described in (\ref{fix-point-iteration-P}); 

\item Plot $\| p_{\varepsilon }^{m} -P^{H_{\varepsilon }^{m}} \|$ versus $t^m$ for each $\varepsilon$ (the norm is taken at a single time layer, chosen as $L^1$ norm in our numerical test). AP property means that after a short period of time (the length of which decreases with $\varepsilon$), the difference between the solution $p_{\varepsilon}^m$ and the corresponding discretized $v$-wise quasi-steady state $P^{H_\varepsilon^m}$ should drop to a small value which vanishes as $\varepsilon \to 0$, i.e. 
\begin{equation} \label{numerical AP convergence}
\| p_{\varepsilon }^{m} -P^{H_{\varepsilon }^{m}} \| \xrightarrow{\varepsilon \rightarrow 0} 0.
\end{equation}
%and AP property of a numerical scheme can be measured by how much (\ref{numerical AP convergence}) is satisfied.
\end{enumerate}
In our simulations, we choose $V_F=2, V_R=1, V_{\min} = -4, a = 1, W_{\min} = -1.1, W_{\max} = 0.1, T_{\max} = 0.3, \sigma(\overline N) = \overline N, \Delta t = 5\times 10^{-4}, \Delta v = 0.1, \Delta w=0.01 $ ,
\begin{equation}\label{p-init}
p_{\text{init}} = \begin{cases}
\displaystyle \sin ^2(\pi v) \sin ^2 (\pi w) & -1 < w < 0 \text{ and } -1 < v < 1\\
0 & \text{otherwise}
\end{cases}
\end{equation}and
\begin{equation}\label{douple_peak_input}
I(w) = \frac{1}{2}\mathrm{e}^{-( 10w+5)^{2}}
\end{equation}

\begin{figure}[htbp]
    \centering
    \includegraphics[scale=0.75]{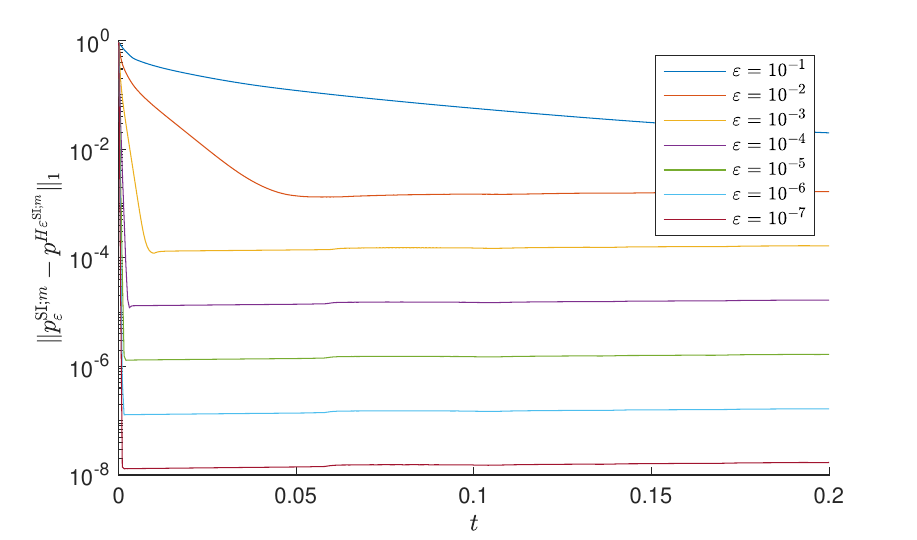}
    \caption{AP property of the $v$-FI scheme. $\| p_{\varepsilon }^{\text{FI};m} -P^{H_{\varepsilon }^{\text{FI};m}}\|$ is plotted over $t \in [0,T]$, discretized with the time step $\Delta t=5\times 10^{-4}$. The discretized quasi-steady state $P^{H_{\varepsilon }^{\text{FI};m}}$ is obtained by the iterations described in (\ref{fix-point-iteration-P}). It can be roughly seen that when $\varepsilon$ is small, $\| p_{\varepsilon }^{\text{FI};m} -P^{H_{\varepsilon }^{\text{FI};m}}\|=o(\varepsilon)$ after a few time steps.}
    \label{fig:AP v-FI}
\end{figure}

\begin{figure}[htbp]
    \centering
    \includegraphics[scale=0.75]{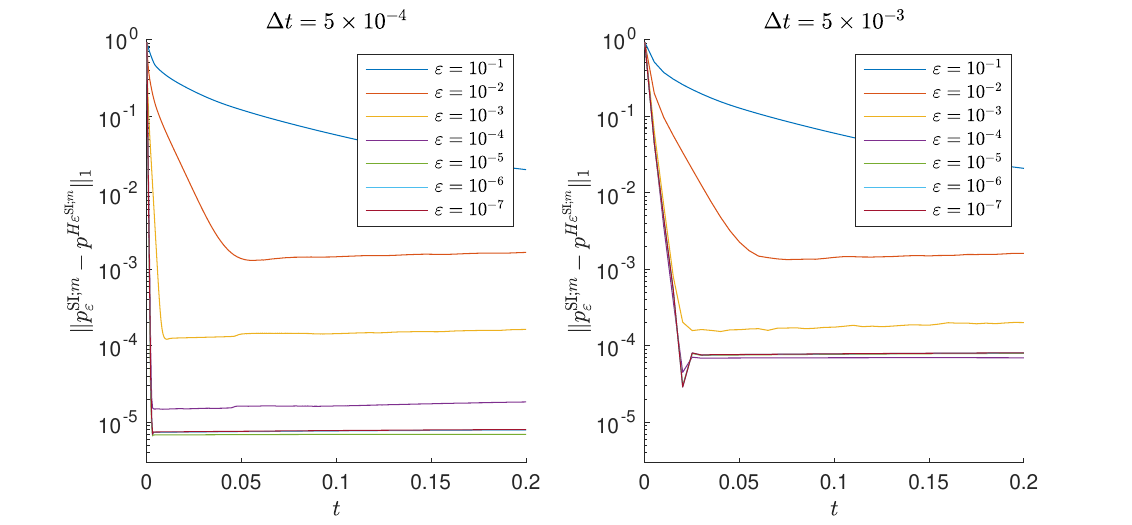}
    \caption{AP property of the $v$-SI scheme. $\| p_{\varepsilon }^{\text{SI};m} -P^{H_{\varepsilon }^{\text{SI};m}}\|$ is plotted over $t \in [0,T]$, discretized with the grid size $\Delta t=5\times 10^{-4}$ (left) and $5 \times 10^{-3}$ (right). The discretized quasi-steady state $P^{H_{\varepsilon }^{\text{SI};m}}$ is obtained by the iterations described in (\ref{fix-point-iteration-P}). Comparing to the $v$-FI scheme, the property $\| p_{\varepsilon }^{\text{SI};m} -P^{H_{\varepsilon }^{\text{SI};m}}\|$ is of order $\varepsilon$ when $\varepsilon$ is comparable to or greater than $\Delta t$. When $\varepsilon \ll \Delta t$ (which is the case when $\varepsilon=10^{-5}, 10^{-6}, 10^{-7}$ in the figure, represented by the green, cyan and blue curves, respectively), $\|p_{\varepsilon }^{\text{SI};m}-P^{H_{\varepsilon }^{\text{SI};m}}\|$ becomes roughly of order of $\Delta t$, so smaller $\varepsilon$ does not make this difference uniformly smaller as it does for the $v$-FI scheme. %\hq{This caption is slightly reorganized to clarify what these numerical results show}
    }
    \label{fig:AP v-SI}
\end{figure}

We perform the aforementioned AP test for both the $v$-FI scheme and the $v$-SI scheme, choosing $\varepsilon$ ranging from $10^{-1}$ to $10^{-7}$, and the results are shown in Figure \ref{fig:AP v-FI} and Figure \ref{fig:AP v-SI}, respectively.

From Figure \ref{fig:AP v-FI}, we can see that for the $v$-FI scheme, the differences between $p_{\varepsilon }^{m}$ and $P^{H_{\varepsilon }^{m}}$ is decreasing when $\varepsilon\to 0$, we can even assert $\| p_{\varepsilon }^{m} -P^{H_{\varepsilon }^{m}} \|=o(\varepsilon)$ from the figure.

From Figure \ref{fig:AP v-SI}, we can see that for the $v$-SI scheme, the differences between $p_{\varepsilon }^{m}$ and $P^{H_{\varepsilon }^{m}}$ is  also  decreasing when $\varepsilon\to 0$. However, this decreasing trend is halting when $\varepsilon = 10^{-5}$ for $\Delta t = 5\times 10^{-4}$ and $\varepsilon = 10^{-4}$ for $\Delta t = 5\times 10^{-3}$, %\jh{Can you explain why they saturate at different $\Delta t$?} \hq{I expanded my explanation in the rest of this section. Is this enough to explain the saturation?} 
and further decrease of $\varepsilon$ does not uniformly draw $p_{\varepsilon }^{m}$ and $P^{H_{\varepsilon }^{m}}$ closer. This can be explained by the time-delayed effect of the $v$-SI scheme, which was presented in section \ref{subsec:asymptotic preserving}. This time-delay effect introduces an $o(\Delta t)$ difference between $p_{\varepsilon }^{m}$ and $P^{H_{\varepsilon }^{m}}$, which will not further decrease when $\varepsilon$ becomes a small number comparing to $\Delta t$. We may say that $\| p_{\varepsilon }^{m} -P^{H_{\varepsilon }^{m}} \|=o(\Delta t)$ when $\varepsilon\ll \Delta t$ and $\| p_{\varepsilon }^{m} -P^{H_{\varepsilon }^{m}} \|=o(\varepsilon)$ when $\varepsilon \gg \Delta t$ for the $v$-SI scheme. %Such a difference from the $v$-wise quasi-steady state may be recognized as a part of discretization error of our scheme.

\subsection{Learning and reacting}
\label{subsec: Learning and Reacting}
From this subsection on, we explore model (\ref{auGi4Tx7zhb}) with our scheme. We always use the $v$-SI scheme in the rest of this section.

In this subsection, we test the model's learning and discrimination abilities. We let the model perform learning-testing tasks that was originally proposed in \cite{BP18}, containing a learning phase and a testing phase:

\noindent\textbf{Learning Phase} 
A heterogeneous input $I(w)$ is presented to the system, while the learning process is active. The learning rule is determined for the inhibitory weights by $-N(w)\overline{N}$  by taking $K(w) = -1$ if $w \leq 0$. After some time, the synaptic weight distribution $H(w,t)$ converges to an equilibrium distribution $H_I^*(w)$, which depends on $I$.

\noindent\textbf{Testing Phase} 
A new input $J(w)$ is presented to the system. This time, the neural system reacts to $J$ but {does not learn it. In this stage, there is only $v$-wise transport and no $w$-wise convection}. To achieve this goal, we implement the $v$-wise quasi-steady state solver (\ref{fix-point-iteration-P}) with $I$ replaced by the testing signal $J$, during which we keep $H_I^*(w)$ we obtained in the Learning Phase.

The pattern recognition property of the system is that it can detect whether the new input $J(w)$ is actually the same one that has been presented during the learning phase, i.e. $I(w)$: indeed, in this case, $N_{I,J}^*(w) = w\boldsymbol 1[-A,0)$ has a very specific shape that does not depend upon the original input $I$ that has been learned in the learning phase.

For the learning phase, We implement the $v$-SI scheme; for the testing phase, we implement the $v$-wise quasi-steady state solver (\ref{fix-point-iteration-P}), replacing the input function $I(w)$ with the testing input function $J(w)$.

In our numerical experiments, the learning and testing input functions are chosen from the set
\begin{equation} \label{input functions}
\{\psi_i(10w+5) + 1 | i = 0, 1, \cdots, 4\}
\end{equation}
where $\psi_{i}(x)$ stands for the $i$-th normalized Hermite function defined recursively as
\begin{equation}
\begin{array}{c}
\displaystyle \psi_{0}(y) \equiv \pi^{-\frac14} \mathrm{e}^{-\frac{1}{2}y^2} ; \quad \psi_{1}(y) \equiv \pi^{-\frac14} \sqrt{2} y \mathrm{e}^{-\frac{1}{2}y^2} \\
\displaystyle \psi_{n+1}(y)=\sqrt{\frac{2}{n+1}} y \psi_{n}(y)-\sqrt{\frac{n}{n+1}} \psi_{n-1}(y).
\end{array}
\end{equation}
This set of functions are well known for forming a set of bases for the Hilbert space $L^2(-\infty, +\infty)$. We take the argument $10w+5$ so that the input functions are centralized at $w=-1/2$ and almost entirely supported in our truncated solving region; and the extra ``$+1$'' term widens the support of the solution so that it can react to input signals at wider range of $w$. We choose $V_F=2, V_R=1, V_{\min} = -4, a = 1, \varepsilon = 0.1,  W_{\min} = -1.1, W_{\max} = 0.1, T_{\max} = 5, \sigma(\overline N) = \overline N, \Delta t = 0.005, \Delta v = 0.1, \Delta w=0.01 $, and $p_{\text{init}}$ is given as (\ref{p-init}). 

\begin{figure}[htbp]
    \centering
    \includegraphics[scale=0.8]{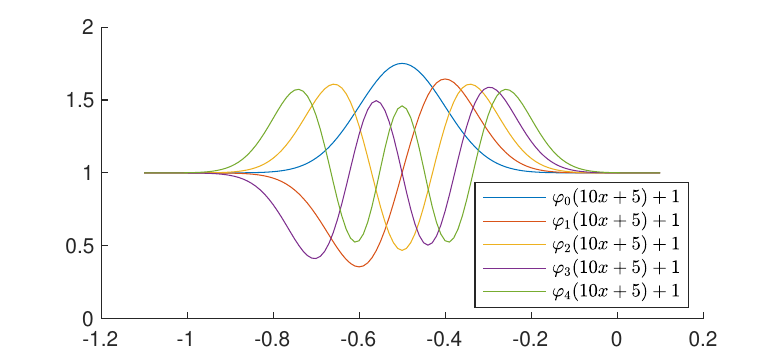}
    \caption{The candidate input functions $\varphi_{i}(10w+5) + 1$ where $i=0,1,2,3,4$ and $\varphi_{i}$ stands for the $i$-th normalized Hermite functions.}
    \label{fig:Hermite Functions}
\end{figure}

\begin{figure}[htbp]
    \centering
    \includegraphics[scale=0.82]{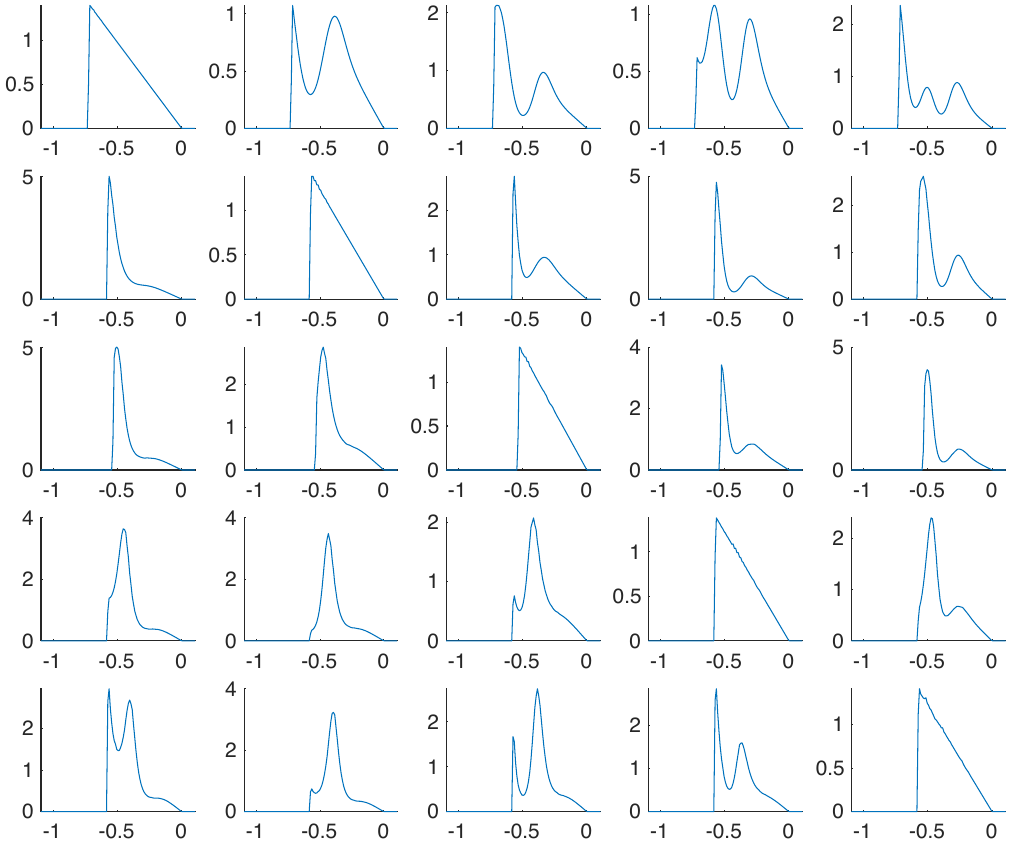}
    \caption{Test on the neural network's recognition property. The pattern recognition task described at the beginning of subsection \ref{subsec: Learning and Reacting} is done for all possible pairs of input functions given in (\ref{input functions}). For each $i,j=0,1,\cdots,4$, the subfigure at the $(i+1)$-th row, $(j+1)$-th column gives the final firing pattern $N_{I,J}^*(w)$ versus after learning the signal $I(w) = \varphi_i(10w+5) + 1$ and reacting to the signal $J(w) = \varphi_j(10w+5) + 1$. For each subfigure, the horizontal coordinate represents $w$, and the vertical axis represents $N(w)$ at the end of the task. It can be seen that, when the learned input $I(w)$ and the testing input $J(w)$ are the same (corresponding to the diagonal entries in the figure), the firing pattern is of the nice triangular shape $N_{I,J}^*(w) = w\boldsymbol 1[-A,0)$; but when $I(w)$ and $J(w)$ are different, the firing pattern will not be regular-shaped.}
    \label{fig: learn-kernel-result}
\end{figure}

For each $i,j=0,1,\cdots,4$, we let the model perform the fore-mentioned task with $I(w) = \psi_i(10w+5)$ and $J(w) = \psi_j(10w+5)$, and the corresponding network activity $N(w)$ is shown in the $(i+1)$-th row, $(j+1)$-th column in Figure \ref{fig: learn-kernel-result}. It can be seen that, when the learned input $I(w)$ and the testing input $J(w)$ are the same (corresponding to the diagonal entries in the figure), the firing pattern is of the nice triangular shape $N_{I,J}^*(w) = w\boldsymbol 1[-A,0)$; but when $I(w)$ and $J(w)$ are different, the firing pattern is not in a regular shape. It is worth noting that in \cite{BP18}, the authors provided a primitive numerical result for the same learning-testing task, the networks firing pattern after reacting to the same learned input signal, shown in the article's Figure 3, turned out to be quite zigzag. In comparison, with our numerical method, an almost-perfect triangular-shaped pattern is produced after the model learn and react to the same input (see the subfigures in the diagonal entries in Figure \ref{fig: learn-kernel-result}). 

\subsection{Positive synaptic weights}
\label{subsec: Positive Synaptic Weights}
This subsection focus on exploring model's behavior when \blue{the neural system is excitatory, i.e. the solution is supported in the with region with $w\geq 0$}. In \cite{BP18}, the authors' analytic work focus mainly on inhibitory networks. And for the excitatory neural system, they provided that, under some choice of parameters, the steady-state solution may not exist. We further explore with our numerical scheme the behavior of the solution when the synaptic weights are positive.

Our numerical exploration shows that in the excitatory cases, \blue{the network} can either develop to a steady state like it did in the inhibitory cases, or exhibit a trend of accelerating expansion unlike it did with inhibitory cases under different parameter settings.

\subsubsection{Long-time steady solution}
\label{subsec: Long-Time Steady Solution}
We firstly test some scenarios in which the solutions present behavior similar to it would present with negative synaptic weights $w$. 

We choose $I\equiv 1$, $a=1$, $\varepsilon=0.2$, $\sigma(\overline N) = \overline N/(1+\overline N)$, $\Delta t = 5\times 10^{-3}, \Delta v = 0.1, \Delta w = 0.01$, $T_{\max} = 5$, $K(w) \equiv 1$ and initial condition given as (\ref{p-init}). The solution finally develops to a steady-state with $N^*(w)$ forming a right triangular pattern supported in a single interval $w \in [0, A]$ for some $A>0$, just like it did with negative synaptic weights. This is the case that corresponds to the correct biological picture.

\begin{figure}[htbp]
    \centering
    \includegraphics[scale=0.9]{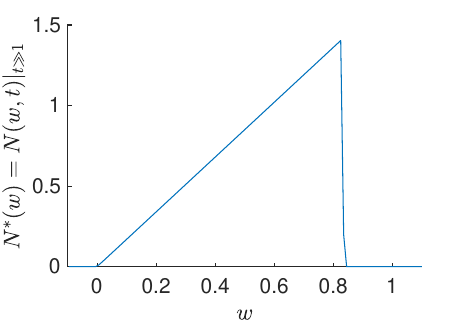}
    \caption{Solution with positive synaptic weight. With $I=0$, $a=1$, $\varepsilon=0.2$, initial condition given as (\ref{p-init}), $\sigma(\overline N) = \overline N/(1+\overline N)$, $\Delta t = 5\times 10^{-3}, \Delta v = 0.1, \Delta w = 0.01$, $T_{\max} = 5$, $K(w) \equiv 1$. The solution finally develops to a steady-state with $N^*(w)$ forming the right triangular pattern just like it did with negative synaptic weights.},
    \label{fig:positive_w_long_time_steady}
\end{figure}

\begin{figure}[htbp]
    \centering
    \includegraphics[scale=0.82]{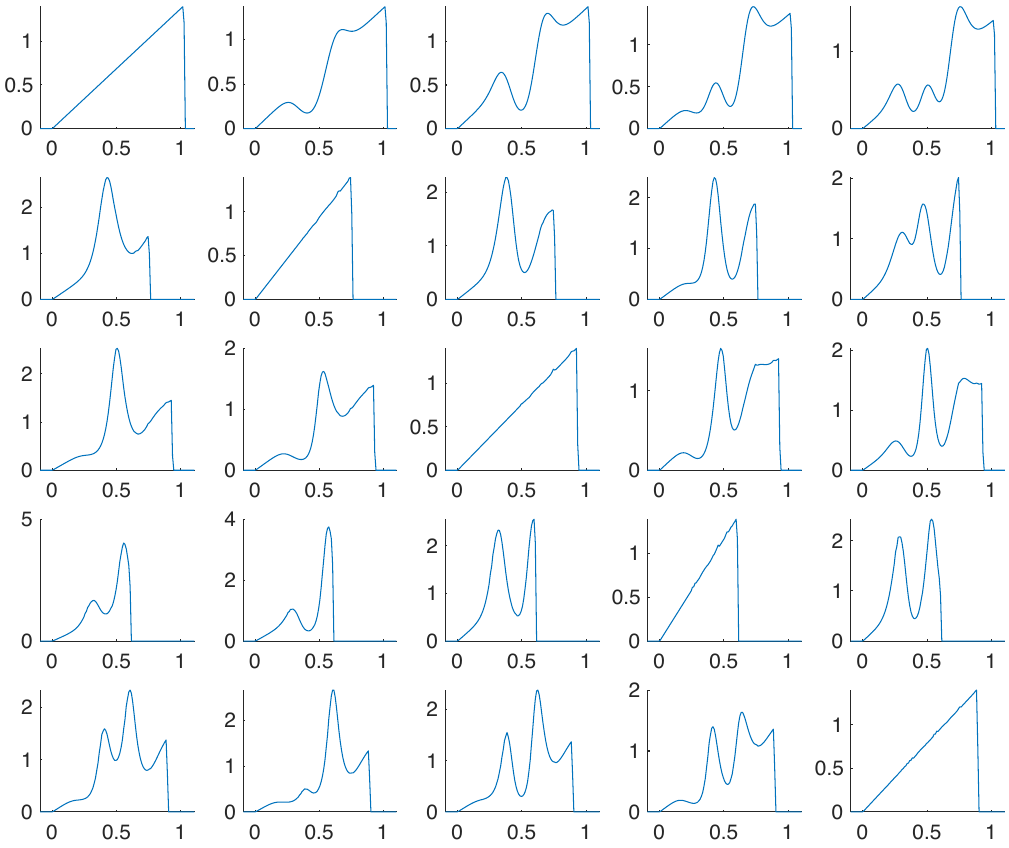}
    \caption{Test of recognition property of excitatory neural networks. The pattern recognition task described at the beginning of subsection \ref{subsec: Learning and Reacting} is done for all possible pairs of input functions given in (\ref{input functions}). For each $i,j=0,1,\cdots,4$, the subfigure at the $(i+1)$-th row, $(j+1)$-th column gives the final firing pattern $N_{I,J}^*(w)$ versus after learning the signal $I(w) = \varphi_i(10w-5) + 1$ and reacting to the signal $J(w) = \varphi_j(10w-5) + 1$. For each subfigure, the horizontal coordinate represents $w$, and the vertical axis represents $N(w)$ at the end of the task. It can be seen that, when the learned input $I(w)$ and the testing input $J(w)$ are the same (corresponding to the diagonal entries in the figure), the firing pattern is of the nice triangular shape $N_{I,J}^*(w) = w\boldsymbol 1(0,A)$; but when $I(w)$ and $J(w)$ are different, the firing pattern will not be regular-shaped.},
    \label{fig:positive_w_learn_test}
\end{figure}

Indeed, as long as it can develop to a steady state, the excitatory neural network also has the recognition ability shown in subsection \ref{subsec: Learning and Reacting}. To show this, we perform a similar learning-testing task as we did in the last subsection. The parameter settings are the same as they were in the last subsection except for that the initial condition $p_{\text{init}}$ is supported in the region with $w\geq 0$, $K(w) \equiv 1$, $\sigma(\overline N) = \overline N / (1+\overline N)$ and the candidate input functions become
\begin{equation}
    \psi_i(10w-5) + 1 \quad i = 0, 1, \cdots, 4
\end{equation}
so that the input are centered at $+1/2$ instead of $-1/2$ for the inhibitory cases. The results are given in Figure \ref{fig:positive_w_learn_test}. The arrangement of the figures are the same as that of Figure \ref{fig: learn-kernel-result}. It can be seen that $w-N(w)$ subfigures at the diagonal entries are in perfect right triangular shape, which is the same as it is in Figure \ref{fig:positive_w_learn_test}.

\subsubsection{Unsteady solution}
\label{subsec: Unsteady Solution}
\begin{figure}
    \centering
    \includegraphics[scale=0.8]{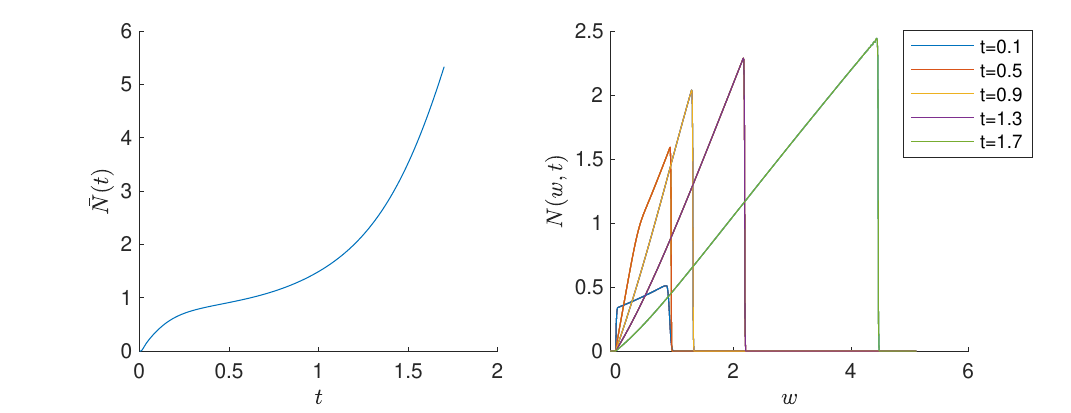}
    \caption{Solution with positive synaptic weight. With $I(w)\equiv 1$, $a=1$, $\varepsilon=0.2$, initial condition given as (\ref{p-init}), $\sigma(\overline N) = 3\overline N/(1+\overline N)$, $\Delta t = 5\times 10^{-3}, \Delta v = 0.1, \Delta w = 0.01$, $T_{\max} = 5$, $K(w) \equiv 1$. Left: in the unsteady scenarios, $\overline N$ grows with accelerate. Right: The solution turns out to propagate to the region with greater $w$. The $w$-wise right boundary of the support of $N$ grows increasingly fast.}
    \label{fig:positive_w_unsteady_N}
\end{figure}
Now we turn our attention to a different case, where the solution is not converging to a steady solution. In \cite{BP18}, it was presented that even if we choose the firing function $\sigma$ to be bounded 
\begin{equation} \label{bounded_sigma}
    \sigma(\overline{N}) = k\frac{\overline N}{1+\overline N},
\end{equation}
we may still face some scenarios that a steady state does not theoretically exist. And we further explore the solution's behavior in this case. 

We choose $I(w) \equiv 1$, $k = 3$, and the rest parameters were same as it was in the first numerical experiment in sub-subsection \ref{subsec: Long-Time Steady Solution}. The solution turns out to propagate to the region with greater $w$. The $w$-wise right boundary of the support of $N$ grows increasingly fast (Figure \ref{fig:positive_w_unsteady_N} right), and henceforth $\overline N$ grows increasingly big (Figure \ref{fig:positive_w_unsteady_N} left). It should be noticed that our result shown in Figure \ref{fig:positive_w_unsteady_N} does not necessarily imply that this is a blow up solution, since we are not sure whether $\overline N(t)$ will grow to $\infty$ in a finite time.

Even if the solution is not a blow-up one, it still does not practically make sense that the network's synaptic weight grows increasingly fast to $\infty$. Therefore, the solution's fast growing behavior is a non-physical effect of our model.

\section{Conclusion}\label{sec:conclusion}
In this work, we have proposed a numerical scheme for the Fokker-Planck equation of structured neural networks (\ref{auGi4Tx7zhb}), which is shown to be conservative, positivity-preserving, and asymptotic-preserving.

The key to success of our scheme is the implicit treatment for flux shift operators in (\ref{auGi4Tx7zhb}), based on which the numerical scheme has an improved positivity preserving property: there is no $\displaystyle v$-wise grid ratio constraint to guarantee the positivity, which is otherwise a major bottleneck for efficient simulation. 

As for numerical exploration, we have conducted the learning-testing numerical experiments and yielded numerical solutions that not only match the theoretical results perfectly, but also fully demonstrate the recognition ability of the model and provide insights for potential uses and future studies. The numerical tests are also extended to the cases of excitatory neural networks. We have also explored with our numerical experiments the behavior of the model when a steady state solution does not exist. These numerical experiments are successful in capturing partial long time trends of the solution but still preliminary. Further numerical experiments can be done to explore, for example, the model's behavior when the input signal $I$ becomes time-dependent, or alternative learning rules that can prevent $\overline{N}(t)$ from growing to infinity for the excitatory cases. 

\section*{Acknowledgement}
J. Hu's research was partially supported by NSF CAREER grant DMS-1654152. Z. Zhou is supported by the National Key R\&D Program of China, Project Number 2020YFA0712000 and NSFC grant No. 11801016, No. 12031013.

\bibliographystyle{unsrt}

\begin{thebibliography}{10}

\bibitem{Synaptic_plasticity_and_memory}
Stephen~J Martin, Paul~D Grimwood, and Richard~GM Morris.
\newblock Synaptic plasticity and memory: an evaluation of the hypothesis.
\newblock {\em Annual review of neuroscience}, 23(1):649--711, 2000.

\bibitem{Hodgkin1952a}
Alan~L Hodgkin and Andrew~F Huxley.
\newblock A quantitative description of membrane current and its application to
  conduction and excitation in nerve.
\newblock {\em The Journal of physiology}, 117(4):500, 1952.

\bibitem{FitzHugh1961Impulses}
Richard FitzHugh.
\newblock Impulses and physiological states in theoretical models of nerve
  membrane.
\newblock {\em Biophysical Journal}, 1(6):445--466, 1961.

\bibitem{lapicque1907recherches}
Louis Lapicque.
\newblock Recherches quantitatives sur l’excitation \'electrique des nerfs
  trait\'ee comme une polarisation.
\newblock {\em J. Physiol. Pathol. Gen}, 9(1):620--635, 1907.

\bibitem{caceres2010analysis}
Mar{\'\i}a~J C{\'a}ceres, Jos{\'e}~A Carrillo, and Beno{\^\i}t Perthame.
\newblock Analysis of nonlinear noisy integrate \& fire neuron models: blow-up
  and steady states.
\newblock {\em The Journal of Mathematical Neuroscience}, 1(1):1--33, 2011.

\bibitem{carrillo2013classical}
Jos{\'e}~A Carrillo, Mar{\'\i}a d~M Gonz{\'a}lez, Maria~P Gualdani, and Maria~E
  Schonbek.
\newblock Classical solutions for a nonlinear fokker-planck equation arising in
  computational neuroscience.
\newblock {\em Communications in Partial Differential Equations},
  38(3):385--409, 2013.

\bibitem{Global-in-time}
María~J. Cáceres, Pierre Roux, Delphine Salort, and Ricarda Schneider.
\newblock Global-in-time classical solutions and qualitative properties for the
  nnlif neuron model with synaptic delay.
\newblock {\em Communications in Partial Differential Equations}, 2018.

\bibitem{brunel2000dynamics}
Nicolas Brunel.
\newblock Dynamics of sparsely connected networks of excitatory and inhibitory
  spiking neurons.
\newblock {\em Journal of computational neuroscience}, 8(3):183--208, 2000.

\bibitem{DGGP17theMeanField}
Grégory Dumont and Pierre Gabriel.
\newblock The mean-field equation of a leaky integrate-and-fire neural network:
  measure solutions and steady states.
\newblock {\em Nonlinearity}, 33, 10 2017.

\bibitem{brunel1999fast}
Nicolas Brunel and Vincent Hakim.
\newblock Fast global oscillations in networks of integrate-and-fire neurons
  with low firing rates.
\newblock {\em Neural computation}, 11(7):1621--1671, 1999.

\bibitem{delarue2015global}
Fran{\c{c}}ois Delarue, James Inglis, Sylvain Rubenthaler, and Etienne
  Tanr{\'e}.
\newblock Global solvability of a networked integrate-and-fire model of
  mckean--vlasov type.
\newblock {\em The Annals of Applied Probability}, 25(4):2096--2133, 2015.

\bibitem{delarue2015particle}
Fran{\c{c}}ois Delarue, James Inglis, Sylvain Rubenthaler, and Etienne
  Tanr{\'e}.
\newblock Particle systems with a singular mean-field self-excitation.
  application to neuronal networks.
\newblock {\em Stochastic Processes and their Applications}, 125(6):2451--2492,
  2015.

\bibitem{liu2020rigorous}
Jian-guo Liu, Ziheng Wang, Yuan Zhang, and Zhennan Zhou.
\newblock Rigorous justification of the fokker-planck equations of neural
  networks based on an iteration perspective.
\newblock {\em arXiv preprint arXiv:2005.08285}, 2020.

\bibitem{zhou2021investigating}
Jian-Guo Liu, Ziheng Wang, Yantong Xie, Yuan Zhang, and Zhennan Zhou.
\newblock Investigating the integrate and fire model as the limit of a random
  discharge model: a stochastic analysis perspective.
\newblock {\em Mathematical Neuroscience and Applications}, 1, 2021.

\bibitem{hebb1949organization}
Donald~Olding Hebb.
\newblock {\em The organization of behavior: A neuropsychological approach}.
\newblock John Wiley \& Sons, 1949.

\bibitem{BP18}
Benoît Perthame, Delphine Salort, and Gilles Wainrib.
\newblock {Distributed synaptic weights in a LIF neural network and learning
  rules}.
\newblock {\em Physica D: Nonlinear Phenomena}, 353-354:20–30, Sep 2017.

\bibitem{JinYanFPLAP11}
Shi Jin and Bokai Yan.
\newblock A class of asymptotic-preserving schemes for the fokker-planck-landau
  equation.
\newblock {\em J. Comput. Physics}, 230:6420--6437, 07 2011.

\bibitem{bessemoulin2012finite}
Marianne Bessemoulin-Chatard and Francis Filbet.
\newblock A finite volume scheme for nonlinear degenerate parabolic equations.
\newblock {\em SIAM Journal on Scientific Computing}, 34(5):B559--B583, 2012.

\bibitem{CCH15}
Jos{\'e}~A Carrillo, Alina Chertock, and Yanghong Huang.
\newblock A finite-volume method for nonlinear nonlocal equations with a
  gradient flow structure.
\newblock {\em Communications in Computational Physics}, 17(1):233--258, 2015.

\bibitem{CC70}
JS~Chang and G~Cooper.
\newblock A practical difference scheme for fokker-planck equations.
\newblock {\em Journal of Computational Physics}, 6(1):1--16, 1970.

\bibitem{PZ18}
Lorenzo Pareschi and Mattia Zanella.
\newblock Structure preserving schemes for nonlinear fokker--planck equations
  and applications.
\newblock {\em Journal of Scientific Computing}, 74(3):1575--1600, 2018.

\bibitem{HZ22}
Jingwei Hu and Xiangxiong Zhang.
\newblock Positivity-preserving and energy-dissipative finite difference
  schemes for the fokker-planck and keller-segel equations.
\newblock {\em arXiv preprint arXiv:2103.16790}, 2021.

\bibitem{boscarino2016high}
Sebastiano Boscarino, Francis Filbet, and Giovanni Russo.
\newblock High order semi-implicit schemes for time dependent partial
  differential equations.
\newblock {\em Journal of Scientific Computing}, 68(3):975--1001, 2016.

\bibitem{liu2018positivity}
Jian-Guo Liu, Li~Wang, and Zhennan Zhou.
\newblock Positivity-preserving and asymptotic preserving method for 2d
  keller-segal equations.
\newblock {\em Mathematics of Computation}, 87(311):1165--1189, 2018.

\bibitem{HH20}
Jingwei Hu and Xiaodong Huang.
\newblock A fully discrete positivity-preserving and energy-dissipative finite
  difference scheme for poisson--nernst--planck equations.
\newblock {\em Numerische Mathematik}, 145(1):77--115, 2020.

\bibitem{BCH20}
Rafael Bailo, Jose~A Carrillo, and Jingwei Hu.
\newblock {Fully discrete positivity-preserving and energy-dissipating schemes
  for aggregation-diffusion equations with a gradient flow structure}.
\newblock {\em Commun. Math. Sci.}, 18:1259--1303, 2020.

\bibitem{zhang2021unified}
Yong Zhang, Yu~Zhao, and Zhennan Zhou.
\newblock A unified structure preserving scheme for a multispecies model with a
  gradient flow structure and nonlocal interactions via singular kernels.
\newblock {\em SIAM Journal on Scientific Computing}, 43(3):B539--B569, 2021.

\bibitem{HS20}
Jingwei Hu and Ruiwen Shu.
\newblock A second-order asymptotic-preserving and positivity-preserving
  exponential runge--kutta method for a class of stiff kinetic equations.
\newblock {\em Multiscale Modeling \& Simulation}, 17(4):1123--1146, 2019.

\bibitem{levenqueNumConserv}
Randall~J. LeVeque.
\newblock {\em Numerical methods for conservation laws (2. ed.).}
\newblock Lectures in mathematics. Birkhäuser, 1992.

\bibitem{JH19}
Jingwei Hu, Jian-Guo Liu, Yantong Xie, and Zhennan Zhou.
\newblock A structure preserving numerical scheme for fokker-planck equations
  of neuron networks: Numerical analysis and exploration.
\newblock {\em Journal of Computational Physics}, 433:110195, 2021.

\bibitem{ple77}
Robert~J Plemmons.
\newblock M-matrix characterizations. i—nonsingular m-matrices.
\newblock {\em Linear Algebra and its Applications}, 18(2):175--188, 1977.

\end{thebibliography}

\end{document}